\documentclass[]{interact}

\usepackage{epstopdf}
\usepackage[caption=false]{subfig}

\usepackage[numbers,sort&compress]{natbib}
\usepackage[ruled,vlined]{algorithm2e}
\usepackage{color}
\bibpunct[, ]{[}{]}{,}{n}{,}{,}
\renewcommand\bibfont{\fontsize{10}{12}\selectfont}

\theoremstyle{plain}

\theoremstyle{definition}

\theoremstyle{remark}

\begin{document}

\title{Parallel Strategies for Direct Multisearch}

\author{
\name{S. Tavares\textsuperscript{a}, C. P.
Br\'as\textsuperscript{b}, A. L.
Cust\'odio\textsuperscript{b}\textsuperscript{*}\thanks{\textsuperscript{*}Corresponding
author. Email: {\tt alcustodio@fct.unl.pt}.}, V.
Duarte\textsuperscript{c}, and P. Medeiros\textsuperscript{c}}
\affil{\textsuperscript{\textbf{a}} FCT NOVA;
\textsuperscript{\textbf{b}}Department of Mathematics, FCT NOVA --
CMA; \textsuperscript{\textbf{c}} Department of Computer Science,
FCT NOVA -- NOVA LINCS.  Campus de Caparica, 2829-516 Caparica,
Portugal.} }

\maketitle

\begin{abstract}

Direct Multisearch (DMS) is a Derivative-free Optimization class
of algorithms suited for computing approximations to the complete
Pareto front of a given Multiobjective Optimization problem. It
has a well-supported convergence analysis and simple
implementations present a good numerical performance, both in
academic test sets and in real applications. Recently, this
numerical performance was improved with the definition of a search
step based on the minimization of quadratic polynomial models,
corresponding to the algorithm BoostDMS.

In this work, we propose and numerically evaluate the performance
of parallelization strategies for this solver, applied to the
search and to the poll steps. The final parallelized version not
only considerably decreases the computational time required for
solving a Multiobjective Optimization problem, but also increases
the quality of the computed approximation to the Pareto front.
Extensive numerical results will be reported in an academic test
set and in a chemical engineering application.

\end{abstract}


\begin{amscode}
90C29, 90C56, 68W10, 90C30.
\end{amscode}

\section{Introduction}\label{sec:introduction}

Direct MultiSearch (DMS)~\cite{ALCustodio_et_al_2011} is a
well-established class of Multiobjective Derivative-free
Optimization algorithms, which has been successfully used in real
applications from very different scientific areas. Optimization of
composite structures~\cite{DAPereira_et_al_2021}, economic
problems~\cite{DFioriti_et_al_2020}, or data science
applications~\cite{ALi_BXue_MZhang_2020} are just a few examples
of the practical value of DMS. Additionally, the corresponding
solver has been considered in benchmark studies, as a reference
code for new
algorithms~\cite{GLiuzzi_SLucidi_FRinaldi_2016,GCocchi_et_al_2018}.

When the algorithmic class of Direct Multisearch was originally
proposed, the goal was to generalize directional direct
search~\cite{CAudet_WHare_2017,ARConn_KScheinberg_LNVicente_2009}
to Multiobjective Optimization. Surprisingly, the first
implementation developed, which did not even comprise a defined
search step, showed to be competitive against other state-of-art
solvers like evolutionary strategies~\cite{KDeb_et_al_2002},
simulated annealing~\cite{SBandyopadhyay_et_al_2008}, or direct
search algorithms~\cite{CAudet_GSavard_WZghal_2008}. Recently,
exceptional good results were also obtained when applying the
solver to derivative-based
optimization~\cite{RAndreani_et_al_2021}.

However, since its first version, in 2010, new releases
corresponded to small improvements of the original implementation,
with no major changes. Exception occurred
in~\cite{CBras_ALCustodio_2020}, where the definition of a search
step based on the minimization of quadratic polynomial models was
proposed. This new version of the code was named as BoostDMS.

This paper considers the use of BoostDMS to solve the
Multiobjective Derivative-free Optimization problem
\begin{equation}\label{eq:MOO}
\begin{split}
\min\quad & F(x)\equiv\left(f_1(x),\dots,f_q(x)\right)^\top\\
\mbox{s.t.} \quad&  lb\leq x\leq ub,
\end{split}
\end{equation}
where $q\geq 2$, $lb<ub$ represent bounds on the problem
variables, and each $f_i:\mathbb{R}^n\rightarrow\mathbb{R}\cup
\{+\infty\},i=1,2,\ldots,q$ denotes a component of the objective
function, for which derivatives are not available, neither can be
numerically approximated.

When the different components of the objective function are
conflicting between each other, it is impossible to find a point
that minimizes every function component. Thus, rather than a
single point, the solution of a Multiobjective Optimization
problem is a set of points, for which it is not possible to
improve the value of one component of the objective function
without deteriorating the value of at least another. Points with
these characteristics are said to be \emph{efficient points} and
the corresponding images form the \emph{Pareto front} of the
Multiobjective Optimization problem.

Since we are addressing Derivative-free Multiobjective
Optimization problems, like in the single objective case, the
computational time required for computing a solution is usually
large, due to the expensive process of function evaluation. In
fact, since the problem solution is no longer a point but a set of
points, there is a strong need for strategies that allow to reduce
the computational time without deteriorating the quality of the
solution computed.

Distributed computing with multiple processors is often considered
in the context of single objective optimization, for solving
computationally quite demanding problems. Parallelization of
Derivative-free Optimization methods for single objective
optimization, particularly direct search ones, has been largely
studied in the literature. Dennis and
Torczon~\cite{JDennis_VTorczon_1991} described a first parallel
version of a pattern search method which evaluates the function in
parallel and synchronizes at each iteration to compare function
values and make updates. A first asynchronous parallel version of
pattern search was proposed
in~\cite{PHough_TGKolda_VTorczon_2001}, dynamically initiating
actions in response to events. Asynchronous approaches were also
followed in~\cite{CAudet_JDennis_SLDigabel_2008}, where processes
solve problems over subsets of variables. Parallel approaches were
also adopted for global single objective direct search
algorithms~\cite{LTWatson_CABaker_2001,JHe_et_al_2009}, often
hybridizing metaheuristics with local direct search
methods~\cite{JDGriffin_TGKolda_2010,AIFVaz_LNVicente_2007}.

In Multiobjective Optimization, parallelization appears typically
associated to metaheuristics, like is the case of evolutionary
strategies~\cite{JBranke_et_al_2004,BCao_et_al_2017,KDeb_PZope_AJain_2003}
or particle swarm
approaches~\cite{SFan_JChang_2009,YZhou_YTan_2011}. A good
overview of different parallel models for Multiobjective
Optimization can be found in~\cite{EGTalbi_et_al_2008}. However,
to the best of our knowledge, parallel approaches have not yet
been considered for Multiobjective Directional Direct Search
algorithms. Thus, the main goal of this paper will be to develop a
parallel version of BoostDMS, which allows a considerable decrease
of the computational time required for solving a problem, without
deteriorating the quality of the computed approximation to the
Pareto front.

The paper is organized as follows. Section~\ref{sec:BoostDMS}
revises the algorithmic structure of BoostDMS algorithm, pointing
out possible steps for parallelization.
Section~\ref{sec:numerical_results} details the different parallel
strategies considered, evaluating the corresponding numerical
performance. Numerical experiments on a chemical engineering
problem are reported in Section~\ref{sec:styrene}, enhancing the
benefits of the selected parallel strategy. Finally, conclusions
are drawn in Section~\ref{sec:conclusions}.

\section{BoostDMS implementation of Direct Multisearch}\label{sec:BoostDMS}

BoostDMS~\cite{CBras_ALCustodio_2020} is a numerical
implementation of a Direct
Multisearch~\cite{ALCustodio_et_al_2011} algorithm, where the
individual or joint minimization of quadratic polynomial models,
built for each component of the objective function, is used in the
definition of a search step.

The algorithm initializes with a list, $L_0$, of feasible,
nondominated points, corresponding function values and stepsize
parameters. In BoostDMS, following the initialization proposed
in~\cite{ALCustodio_et_al_2011}, $n$ points are evenly spaced in a
line joining the problem bounds. Additionally, the centroid of the
feasible region is always considered as initialization. This list
represents the current approximation to the Pareto front of the
problem and is updated every time that a new feasible nondominated
point is found, corresponding to a successful iteration. The new
nondominated point is added to the list and all the dominated
points are removed from it. Using the strict partial order induced
by the cone $\mathbb{R}_{+}^q$, a point $x$ is said to dominate a
point $y$ if $F(x)\prec_F F(y)$, i.e., when $F(y)-F(x) \in
\mathbb{R}_{+}^q\setminus\{0\}$.

Each iteration of the algorithm starts by selecting an iterate
point, $(x_k, F(x_k),\alpha_k)$, from the current list of feasible
nondominated points, $L_k$, corresponding to the largest gap in
the current approximation to the Pareto front, measured using the
$\Gamma$ spread metric: \begin{equation}\label{gamma_metric}
 \Gamma \; = \; \max_{i \in
\{1,\dots,q\}}\left(\max_{j\in \{1,\dots,N-1
\}}\{\delta_{i,j}\}\right).
\end{equation}
\noindent Here $\delta_{i,j}= f_{i}(y_{j+1})-f_{i}(y_j)$, assuming
that for each component $i$ of the objective function,
$y_1,\ldots,y_{N}$ are the current points in the list $L_k$,
sorted by increasing order of the corresponding objective function
component value. Ties are broken by considering the point with the
largest stepsize parameter. This iterate point will be used at the
search step, as center of the quadratic polynomial models that
will be built for each component, $f_i$, of the objective
function. In case of failure of the search step in computing a new
feasible nondominated point, the poll step will be performed, with
the iterate point as poll center.

The search step reuses previously feasible evaluated points (not
necessarily nondominated), kept in a cache, thus with no
additional cost in terms of objective function values required for
model computation. Points are selected in a neighborhood of the
iterate point $x_k$, i.e., in $B(x_k;\Delta_k)$, where $\Delta_k$
is proportional to the stepsize parameter $\alpha_k$. Models are
built once that $n+2$ points have been evaluated in this region.
Depending on the number of points available, minimum Frobenius
norm models, determined interpolation models or regression
approaches can be
considered~\cite{ARConn_KScheinberg_LNVicente_2009}.

Let $m_i$ define the quadratic polynomial model centered at $x_k$,
corresponding to the objective function component $f_i$,
\begin{equation} \label{eq:polynomial model}
m_i(x) = f_i(x_k) + {g^i_k}^\top (x-x_k)+\frac{1}{2}(x-x_k)^\top
H^i_k (x-x_k), \, \, i=1,2,\ldots,q
\end{equation}
where the gradient vector $g^i_k$ and the symmetric Hessian matrix
$H^i_k$ are computed by solving the linear system corresponding to
the interpolation conditions
\begin{equation}\label{eq:interp_system}m_i(y_j) = f_i(y_j), \, j=1,\ldots,p,\end{equation}
with $y_j$ representing a feasible point for which the objective
function has been previously evaluated, and $p$ denoting the
number of points inside $B(x_k;\Delta_k)$ that can be used in the
computation.

Models can be minimized individually inside $B(x_k;\Delta_k)$, in
an attempt to improve the ability of the algorithm in generating
approximations to the extreme points of the Pareto front. However,
the joint minimization of models is also considered, using a
weighted Chebyshev norm scalarization, by solving the following
problem:

\begin{equation} \label{eq:MOOsearch}
\begin{array}{rl}
\min & \zeta\\[1ex]
{\rm s.t.} & m_i(x) \; \le \; \zeta, \quad i \in I\\
           & \|x-x_k \| \; \le \; \Delta_k \\
           &  lb\leq x\leq ub,
\end{array}
\end{equation}

\noindent where $I \subseteq\{1,2,\ldots,q\}$. For each
cardinality $1\leq l\leq q$ of a subset $I$ there are $C^q_{l}$
possibilities to jointly minimize the models corresponding to the
different components of the objective function. The algorithm
starts by the individual minimization of each model ($l=1$). If it
fails in finding a new feasible nondominated point for the
objective function, combinations of two models are considered
($l=2$). The process is iteratively repeated, increasing the level
$l$ of model combinations until the joint minimization of all
models is considered ($l=q$). If no new feasible nondominated
point is found for the objective function, the poll step will be
performed.

The poll step corresponds to a local search around the selected
iterate point (the poll center), testing a set of directions with
an adequate geometry, scaled by the stepsize parameter. Typically,
positive spanning sets are considered~\cite{CDavis_1954} that
conform to the geometry of the feasible region around the poll
center. For bound constraints, BoostDMS considers coordinate
directions, which are evaluated in a complete polling approach.

At the end of each iteration the stepsize is updated, keeping it
constant for successful iterations and halving the stepsize of the
poll center for unsuccessful ones.

A simplified description of BoostDMS is provided in
Algorithm~\ref{alg:Boostdms}. For a detailed description
see~\cite{CBras_ALCustodio_2020,ALCustodio_et_al_2011}.

\begin{algorithm}
\begin{rm}
\begin{description}
{\small
\item[Initialization] \ \\
Choose a set of feasible points $\{lb \leq x_{ini}^j\le ub,\, j\in
\Upsilon\}$ with $f_i(x_{ini}^j)<+\infty,\forall i
\in\{1,\ldots,q\},\forall\, j\in \Upsilon$, and
$\alpha_{ini}^j>0,\, j\in \Upsilon$ initial stepsizes. Let
$\mathcal{D}$ be a set of positive spanning sets. Initialize the
cache of previously feasible evaluated points and corresponding
function values
$$L_{cache}=\{ (x_{ini}^j,F(x_{ini}^j)), j\in
\Upsilon\}.$$ Retrieve all nondominated points from $L_{cache}$
and initialize the list of feasible nondominated points,
corresponding function values and stepsize parameters
$$L_0=\{ (x_{ini}^j,F(x_{ini}^j),\alpha_{ini}^j ), j\in
\bar{\Upsilon}\subseteq\Upsilon\}.$$
\item[For $k=0,1,2,\ldots$] \ \\
\begin{enumerate}
\item[1.] {\bf Selection of an iterate point:} Order the list
$L_k$ according to the largest gap measured with the $\Gamma$
metric and select the first item $(x, F(x),\alpha) \in L_k$ as the
current iterate, function value, and stepsize parameter (thus
setting $(x_k,F(x_k),\alpha_k)=(x,F(x),\alpha)$). \item[2.] {\bf
Search step:}\ \\ Select a subset of points in $L_{cache}\bigcap
B(x_k;\Delta_k)$ to build the quadratic \noindent polynomial
models. If the cardinality of the subset is smaller than $n+2$ go
to the poll
step.\ \\
\begin{description}
\item[For $i=1,2,\ldots,q$] \ \\
Build the quadratic polynomial model $m_i$, corresponding to the
objective function component $f_i$.
\item[Endfor]\ \\\hspace{-0.7cm}Set $l=0$.\ \\
\item[While $l < m$] \ \\
Set $l:=l+1$, define $J$ the set of all combinations of $l$
quadratic polynomial models taken from the total set of $q$ models
and set $S= \emptyset$.
\begin{description}
\item[For $j=1,2,\ldots,|J|$] \ \\
Compute the point $s_j$, solution of problem (\ref{eq:MOOsearch})
considering $I$ as the set composed by the polynomial models
corresponding to combination $j$. Update $S = S \cup \{s_j\}$.
\item[Endfor]\ \\
\item[Check for success] \ \\
Evaluate $F$ at each point in $S$ and update $L_{cache}$. Compute
$L_{trial}$ by removing all dominated points from $L_{k} \cup
\{(x_s,F(x_s),\alpha_k): x_s\in S \}$. If $L_{trial} \neq L_k$,
set $L_{k+1}=L_{trial}$, stop the cycle loop \textbf{while},
declare
the iteration as successful and skip the poll step.\ \\
\end{description}
\item[Endwhile]\ \\
\end{description}
\item[3.] {\bf Poll step:} Choose a positive spanning set~$D_k$
from the set $\mathcal{D}$. Evaluate $F$ at the feasible poll
points belonging to $\{ x_k+\alpha_k d: \, d \in D_k \}$ and
update $L_{cache}$. Compute $L_{trial}$ by removing all dominated
points from $L_{k} \cup \{(x_k+\alpha_k d, F(x_k+\alpha_k
d),\alpha_k): d\in D_k \wedge lb\leq x_k+\alpha_k d \leq ub\}$. If
$L_{trial} \neq L_k$, declare the iteration as successful and set
$L_{k+1}=L_{trial}$. Otherwise, declare the iteration unsuccessful
and set $L_{k+1} = L_k$. \vspace{1ex} \item[4.] {\bf Stepsize
parameter update:} If the iteration was unsuccessful halve the
stepsize parameter corresponding to the poll center, replacing
$(x_k,F(x_k),\alpha_k)\in L_{k+1}$ by
$(x_k,F(x_k),\frac{\alpha_k}{2})$.\ \\ \vspace{1ex}
\end{enumerate}
\item[Endfor]\ \\}
\end{description}
\end{rm}
\caption{A simplified description of BoostDMS.}
\label{alg:Boostdms}
\end{algorithm}

The convergence analysis of Direct Multisearch algorithmic class
relies on the behavior of the method at unsuccessful poll
steps~\cite{ALCustodio_et_al_2011}, holding independently of the
initialization, of the strategy considered for the selection of
the iterate point, of the definition of a search step, or of the
type of polling strategy considered (opportunistic variants are
allowed, in which polling is interrupted once that a new feasible
nondominated point is found). Recently, worst-case complexity
bounds were also derived for some variants of Direct
Multisearch~\cite{ALCustodio_et_al_2021}. BoostDMS, as an
algorithmic instance of DMS, inherits the convergence properties
of this general class of Multiobjective Derivative-free
Optimization methods.

From the algorithmic description provided, it is clear that
parallelization strategies can be applied to both the search and
the poll steps, enhancing the numerical performance of the
algorithm. Section~\ref{sec:numerical_results} will detail
possible strategies for parallelization, from the simple
distribution of objective function evaluations among the different
processors to more sophisticated schemes, including the parallel
computation of the quadratic polynomial models, or the
simultaneous selection of several iterate points.

\section{Parallelization strategies}\label{sec:numerical_results}

BoostDMS is part of a Derivative-free Optimization Toolbox of
solvers, suited for local and global single objective and
multiobjective problems, implemented in Matlab and freely
available for use
at~\texttt{https://docentes.fct.unl.pt/algb/pages/boostdfo}, under
a GNU Lesser General Public License. The toolbox provides a GUI,
which allows a non-expert user to take advantage of the different
solvers options, including parallelization~\cite{STavares_2020}.
For running BoostDMS, the Matlab toolboxes of Optimization and
Parallel Computing are required for the minimization of the
quadratic models at the search step and for distributing tasks to
the different processors, respectively.

The expensive cost of function evaluation, associated to
Derivative-free Optimization problems, naturally motivates its
parallelization. However, additional parallelization strategies
can and will be considered, not only with the goal of time
reduction but also with the concern of keeping or improving the
quality of the computed approximation to the Pareto front.

To assess this quality, typical metrics from the Multiobjective
Optimization literature are used: purity and spread, as defined
in~\cite{ALCustodio_et_al_2011}, and the hypervolume
indicator~\cite{EZitzler_LThiele_1998,EZitzler_et_al_2003}. For
completeness, we provide a short description of each one of these
metrics. A recent revision of Multiobjective Optimization metrics
can be found in~\cite{CAudet_et_al_2020}.

Purity measures the capability of a given solver in generating
nondominated points. Let $F_{p,s}$ denote the approximation to the
Pareto front computed for problem $p\in \mathcal{P}$ by solver
$s\in \mathcal{S}$, and $F_p$ denote the true Pareto front of
problem~$p$. The purity metric, for problem $p\in \mathcal{P}$ and
solver $s\in \mathcal{S}$, is defined by the ratio
$$Pur_{p,s} = \frac{| F_{p,s} \cap F_p |}{ | F_{p,s} |},$$
taking values between zero and one. Higher values of $Pur_{p,s}$
indicate a better Pareto front in terms of percentage of
nondominated points generated. In general, $F_p$ is not known,
being considered an approximation to it by joining all the final
approximations computed by all the solvers tested and removing all
the dominated points from it.

Spread metrics assess the quality of the distribution of
nondominated points across the final approximation to the Pareto
front generated by the solver. Since we are interested in
computing the complete Pareto front, spread metrics need to
consider the `extreme points' associated to each objective
function component (see~\cite{ALCustodio_et_al_2011}). The
$\Gamma$ metric~(\ref{gamma_metric}), previously used to select
the iterate point based on the maximum gap between consecutive
points lying on the approximation to the Pareto front, is now
computed using a set that includes the `extreme points'. Let us
assume that solver $s \in {\cal S}$ has computed, for problem $p
\in {\cal P}$, an approximated Pareto front with points
$y_1,y_2,\ldots,y_N$, to which we add the `extreme points'
mentioned above ($y_0$ and $y_{N+1})$. The metric $\Gamma_{p,s} >
0$, for
 problem $p\in \mathcal{P}$ and solver
 $s\in \mathcal{S}$, is given by
\begin{equation}\label{gamma_metric2}
 \Gamma_{p,s} \; = \; \max_{i \in
\{1,\dots,q\}}\left(\max_{j\in
\{0,\dots,N\}}\{\delta_{i,j}\}\right),
\end{equation}
\noindent where $\delta_{i,j}=f_{i}(y_{j+1})-f_{i}(y_j)$ (assuming
that the objective function values have been sorted by increasing
order for each objective function component~$i$).

The second spread metric~\cite{KDeb_et_al_2002} measures the
uniformity of the gap distribution across the computed Pareto
front: \begin{equation}\label{delta_metric}\displaystyle
\Delta_{p,s} \; = \;
\max_{i\in\{1,\dots,q\}}\left(\frac{\delta_{i,0}+\delta_{i,N}+\sum_{j=1}^{N-1}|\delta_{i,j}-
\bar{\delta}_i|}{\delta_{i,0}+\delta_{i,N}+(N-1)\bar{\delta}_i}\right),
\end{equation}
where $\bar{\delta}_i$, for $i=1,\ldots,q$, represents the average
of the distances $\delta_{i,j}$, $j=1,\dots,N-1$.

The hypervolume indicator is a good compromise between spread and
purity, since it measures the volume of the portion of the
objective function space that is dominated by the computed
approximation to the Pareto front and is limited by an upper
corner $U_p \in \mathbb{R}^q$ (a point that is dominated by all
points belonging to the different approximations computed for the
Pareto front by all solvers tested). Thus, the hypervolume
indicator can be defined as follows:
$$HV_{p,s} =  Vol\{y \in \mathbb{R}^q| \, y \le U_p \wedge \exists x \in F_{p,s} : x \le y\} = Vol \left(\bigcup_{x \in F_{p,s}} [x, U_p]\right),$$
where $Vol(.)$ denotes the Lebesgue measure of a $q$-dimensional
set of points and $[x, U_p]$ denotes the interval box with lower
corner $x$ and upper corner $U_p$. The approach proposed
in~\cite{CMFonseca_et_al_2006} was used for its computation, after
scaling all the hypervolume values to the interval $[0, 1]$
(see~\cite{CBras_ALCustodio_2020} for details).

The final metric considered is the CPU time required for solving a
problem. As test set, we considered 99 bound constrained
Multiobjective Optimization problems from the collection available
at \texttt{http://www.mat.uc.pt/dms}, now coded in Matlab, with a
number of variables, $n$, between $1$ and $30$, and with $2, 3$,
or $4$ objective function components (problem \texttt{WFG1} was
not considered due to some numerical instabilities of Matlab).
Since these are academic problems, the time associated to function
evaluation is quite small, of the order of milliseconds, very
different from real applications, where function evaluation
typically requires seconds or even minutes of computational time.
For a fair assessment of the benefits of parallelization
strategies, each function was modified with an average
computational time of 0.1 seconds (in single objective
optimization, computational times of similar magnitude have
already allowed benefits from
parallelization~\cite{STavares_2020}). Results are reported for
average CPU time, considering five runs.

Performance profiles~\cite{EDDolan_JJMore_2002} are computed for
the different metrics. Let $t_{p,s}$ be the value of a given
metric obtained for problem~$p \in \mathcal{P}$ with solver~$s \in
\mathcal{S}$, assuming that lower values of $t_{p,s}>0$ indicate
better performance. The cumulative distribution function for
solver $s \in \mathcal{S}$ is given by:
$$\rho_s(\tau)=\frac{1}{|\mathcal{P}|}|\{p\in
\mathcal{P}:r_{p,s}\leq\tau\}|,$$ with
$r_{p,s}=t_{p,s}/\min\{t_{p,\bar{s}}:\bar{s}\in \mathcal{S}\}$.
Thus, the value of $\rho_s(1)$ represents the probability of
solver~$s$ winning over the remaining ones. Since for purity and
hypervolume larger values indicate better performance, for these
metrics the profiles considered $t_{p,s}:=1/t_{p,s}$, as proposed
in~\cite{ALCustodio_et_al_2011}.

All numerical experiments were performed using Matlab R2020a, on a
Xeon Platinium 8171M CPU machine from the Azure cloud, running
Linux as operating system. Default options were assumed for
BoostDMS~\cite{CBras_ALCustodio_2020} and all executions ran using
8 processors. In particular, a maximum of 20\,000 function
evaluations or a minimum value of $10^{-3}$ for the stepsize
parameter $\alpha_k$ of each point in the list were used as
stopping criteria.

\subsection{Parallelization of objective function
evaluations}\label{sec: Parallel feval}

Typical Derivative-free Optimization problems are associated to
expensive function evaluation. In BoostDMS, parallelization of
function evaluations was considered both at the poll and the
search steps, as a way of improving the numerical performance of
the solver.

In Directional Direct Search, the poll step has a natural
structure for parallelization. BoostDMS considers a complete
polling approach, allowing to take full advantage of an
embarrassingly parallel scheme. The set of points to be evaluated
at the poll step is generated and, before function evaluation,
infeasible points or points considered identical (within a
tolerance equal to the minimum stepsize allowed) to others already
included in the cache are discarded. This filtered set of points
is now evaluated in parallel.

All evaluated points are included in the cache ($L_{cache}$) and
dominance is checked, adding nondominated points to the list
$L_k$, one by one, in a similar way to what is done in the
sequential version. This approach ensures that the approximation
generated for the Pareto front is exactly the same than the one
generated by the sequential version, keeping the quality of the
final computed solution.

As described in Section~\ref{sec:BoostDMS}, the search step
considers an approach by levels, when minimizing the quadratic
polynomial models. Subproblem~(\ref{eq:MOOsearch}) may need to be
solved for all combinations of $l$ models. At level $l=1$, the
individual minimization of the models is achieved through function
\texttt{trust.m}, from Matlab, and projecting the resulting points
in the feasible region. For levels two and above ($l\geq2$), the
joint minimization of models resources to Matlab's function
\texttt{fmincon.m}. All solutions of the subproblems at a given
level are evaluated and, only in case of failure in finding a new
feasible nondominated point, the next level is considered.

Three different approaches have been tested, regarding the
parallelization of function evaluations at the search step. The
first, \emph{Par evals within levels}, keeps the function
evaluation structure by levels of the sequential implementation. A
maximum of $q$ batches of parallel function evaluations can be
performed (one associated to each level). In this case, the
quality of the solution is identical to the one of the sequential
version, although a reduction in computational time is expected.
The second variant, denoted by \emph{Par evals 2 batches},
considers function evaluation in two batches, one corresponding to
level $l=1$ and, in case of failure in finding a new feasible
nondominated point, all the remaining levels, corresponding to the
joint minimization of two or more models. Finally, the variant
denoted by \emph{Par evals no levels} solves all
subproblems~(\ref{eq:MOOsearch}), evaluating the objective
function at the corresponding solutions in parallel. In any of the
three strategies, subproblems~(\ref{eq:MOOsearch}) are always
solved sequentially. Only function evaluations are parallelized,
being the procedure followed identical to the one adopted for the
parallelization of the poll step.

A level-by-level approach would be beneficial if success is often
found in lower levels, even though it could not take advantage of
all the potential of parallelizing function evaluations. However,
if success is typically achieved in higher levels or if the search
step is often unsuccessful, considering a single batch of parallel
function evaluations of points, corresponding to all levels, would
be more efficient.

\begin{figure}[htbp]
\begin{center}
\includegraphics[width=7.0cm,height=5cm]{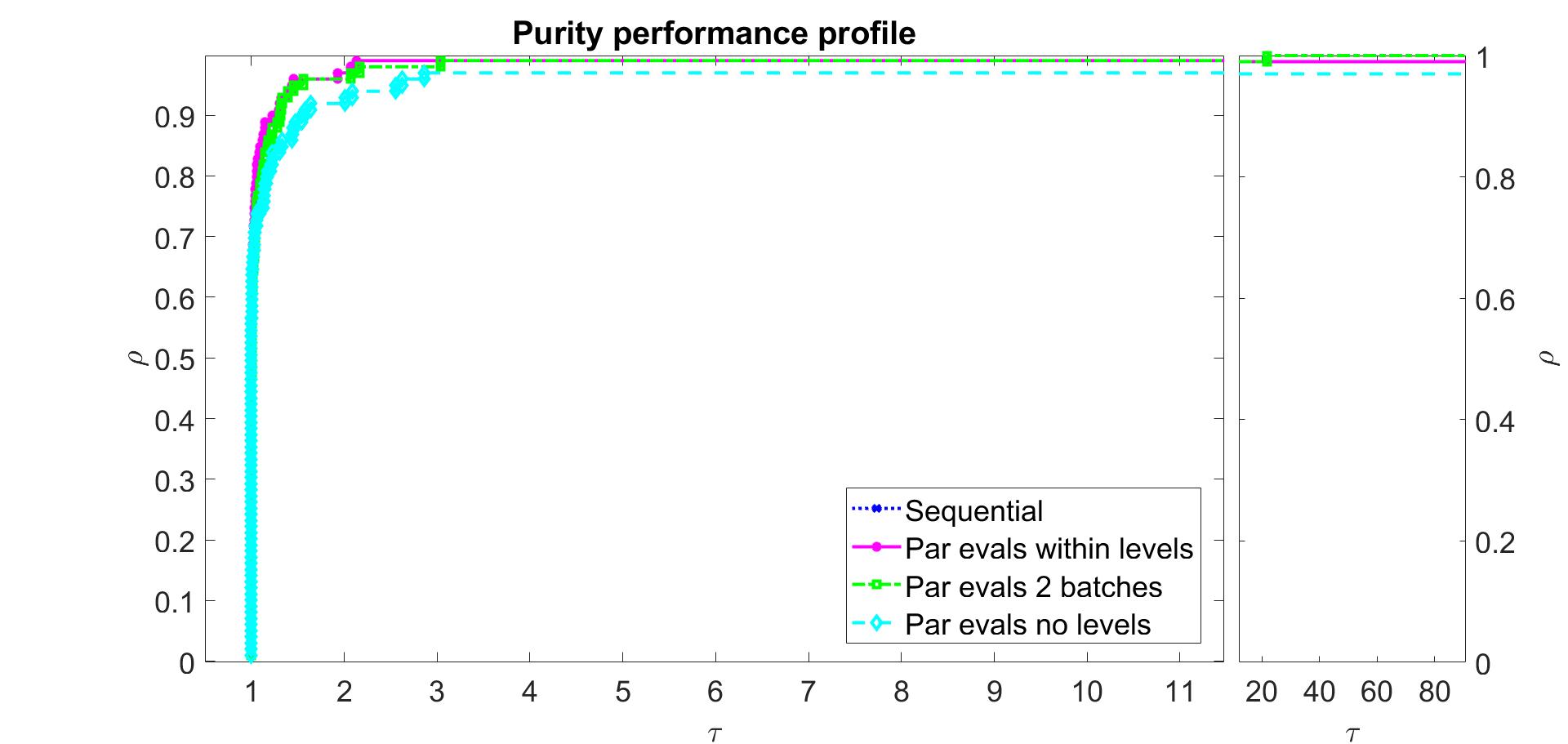}
\includegraphics[width=7.0cm,height=5cm]{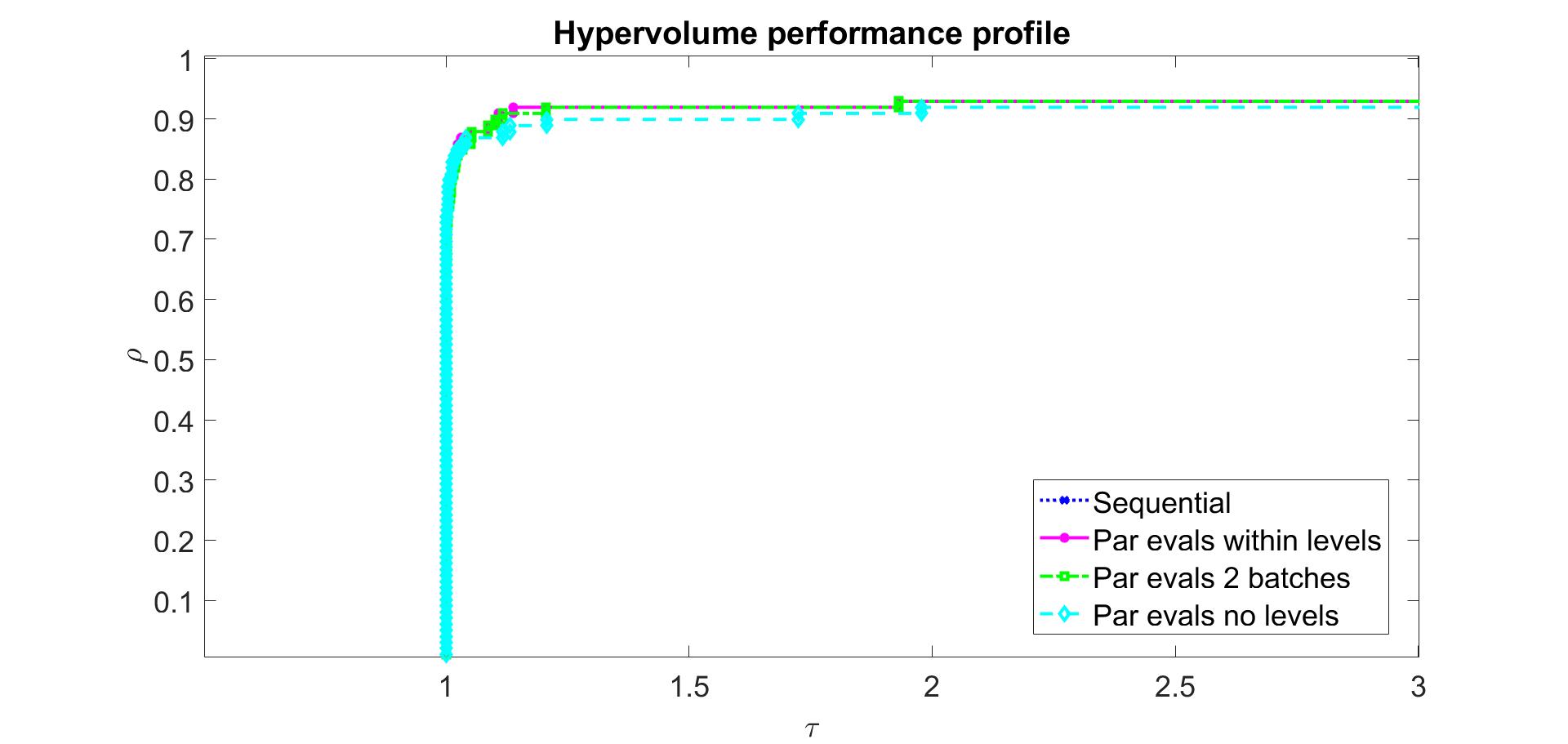}
\end{center}
\caption{\label{fig_seq127_purity_hyper} Comparison between
sequential version of BoostDMS, and parallel versions \emph{Par
evals within levels}, \emph{Par evals 2 batches}, and \emph{Par
evals no levels}, by means of purity and hypervolume metrics.}
\end{figure}

\begin{figure}[htbp]
\begin{center}
\includegraphics[width=7.0cm,height=5cm]{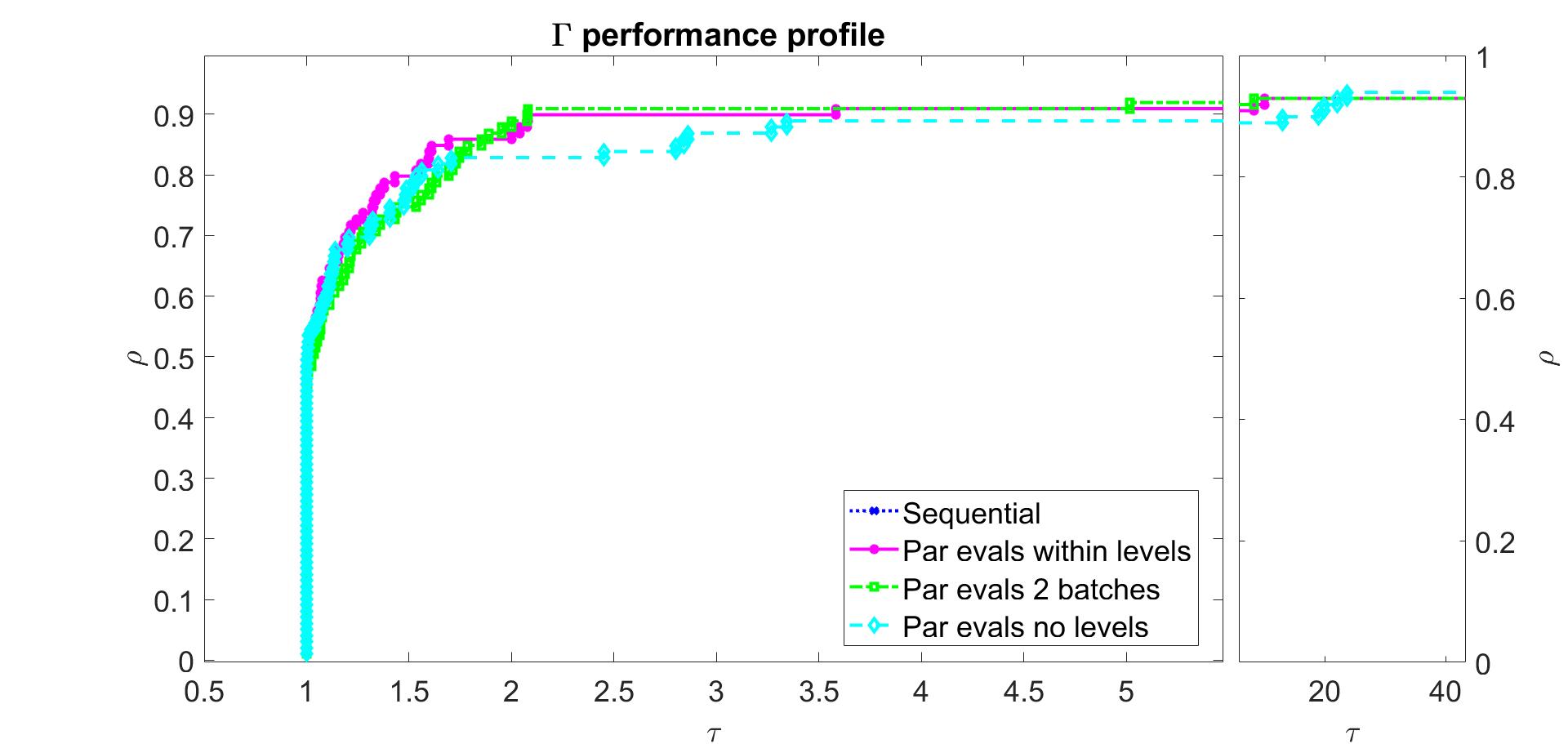}
\includegraphics[width=7.0cm,height=5cm]{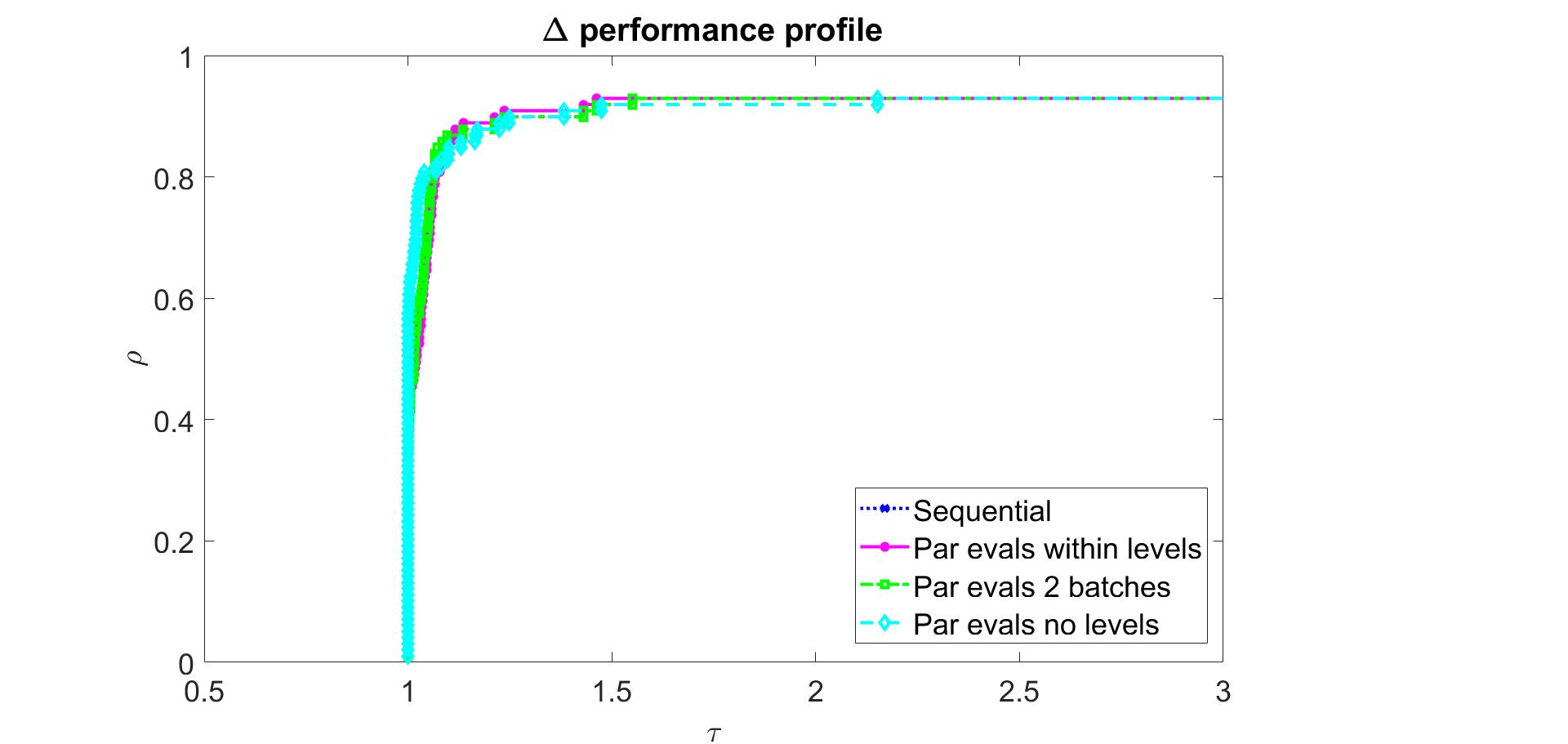}
\end{center}
\caption{\label{fig_seq127_gamma_delta} Comparison between
sequential version of BoostDMS, and parallel versions \emph{Par
evals within levels}, \emph{Par evals 2 batches}, and \emph{Par
evals no levels}, by means of spread metrics ($\Gamma$ and
$\Delta$).}
\end{figure}

\begin{figure}[htbp]
\begin{center}
\includegraphics[width=7.0cm,height=5cm]{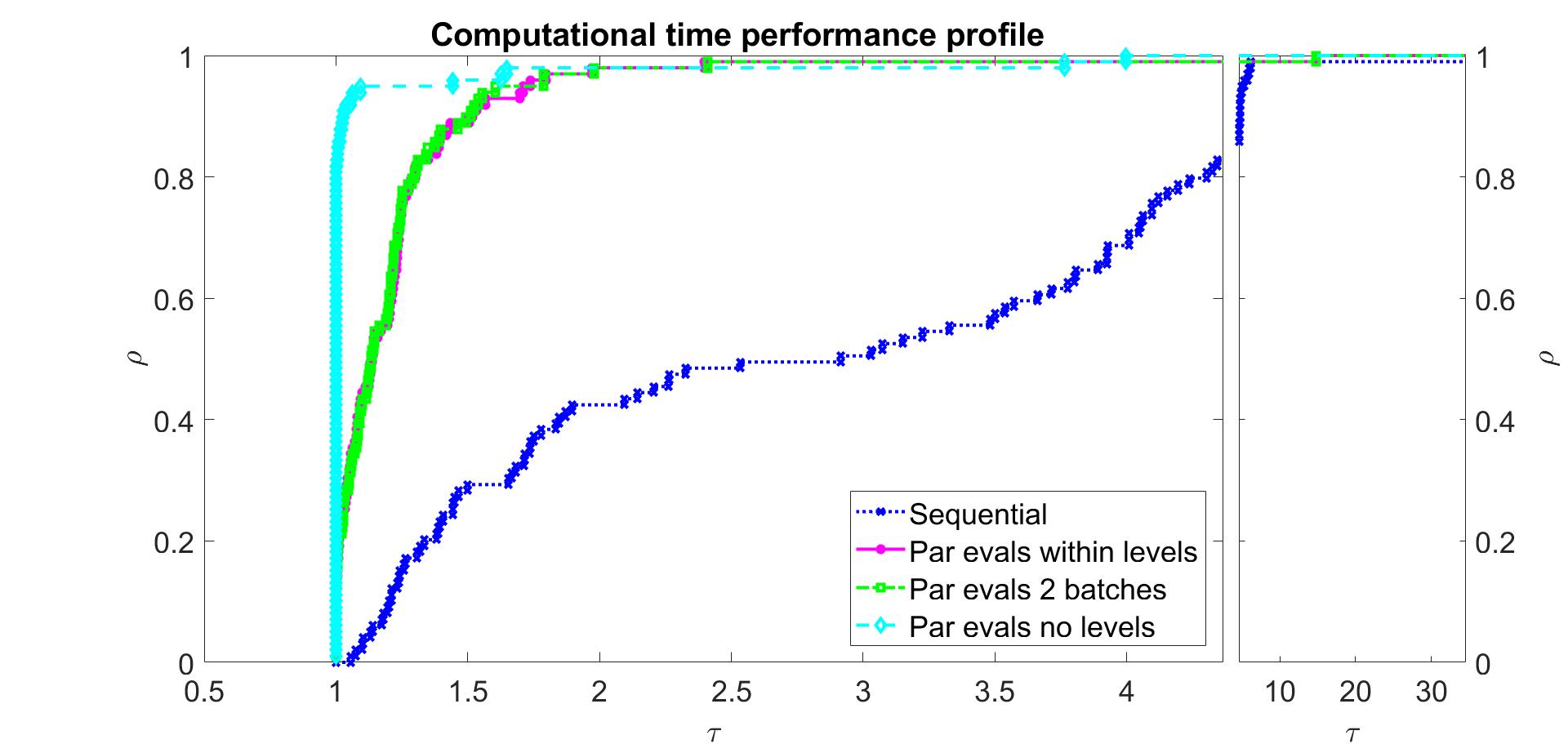}
\end{center}
\caption{\label{fig_seq127_time} Comparison between sequential
version of BoostDMS, and parallel versions \emph{Par evals within
levels}, \emph{Par evals 2 batches}, and \emph{Par evals no
levels}, considering the computational time.}
\end{figure}

Figures~\ref{fig_seq127_purity_hyper},
\ref{fig_seq127_gamma_delta}, and \ref{fig_seq127_time} compare
the results obtained with the sequential version to the three
parallelization strategies described above, for the five metrics
considered.

All strategies present a similar performance in terms of
hypervolume and $\Delta$ metrics. As expected, there is a clear
advantage in computational time for the parallel versions, with
the variant \emph{Par evals no levels} performing the best (see
Figure~\ref{fig_seq127_time}). A strategy that gathers all the
points before evaluation, reduces the number of batches to be
parallelized, translating to better computational times.

However, analyzing the plots corresponding to purity and $\Gamma$
metrics, we conclude that this strategy presents a worse
performance than the other parallel versions. In fact, the results
for the three remaining strategies are very close, indicating that
feasible nondominated points are typically found in the first
level, when the individual minimization of models occurs. Strategy
\emph{Par evals no levels} performs a larger number of function
evaluations per search step, promoting the exhaustion of the
function evaluations budget at early iterations, which could
contribute to better computational times, but does not allow to
refine the quality of the computed approximation to the Pareto
front.

The variant \emph{Par evals within levels}, that keeps the
structure of the parallel function evaluations in batches
corresponding to levels, allows a good performance in terms of
computational time, keeping the quality of the final solution
generated by the sequential version. This was the option taken in
terms of strategy for function evaluation and will be adopted in
the following numerical experiments.

\subsection{Parallelization of models computation and minimization}

The parallelization can obviously be extended to the computation
of the quadratic models for each component of the objective
function and to the corresponding individual or joint
minimization. Thus, three new parallel strategies were developed.

The first, denoted by \emph{Model min (level 1) and par evals},
parallelizes the individual computation and minimization of
models, for each objective function component. Strategy
\emph{Model min (level} $\ge$ \emph{2) and par evals} builds
models and performs their individual minimization sequentially,
only parallelizing the joint minimization of models. Considering
the results of Section~\ref{sec: Parallel feval}, the approach of
evaluating batches of points by level is always followed. Since
the numerical results of Section~\ref{sec: Parallel feval}
indicate that the search step is mainly successful at level $l=1$,
the number of search steps where levels $l\geq 2$ will be
considered is reduced, thus no major improvement is expected with
this second variant. A final strategy, denoted by \emph{Model min
(all) and par evals}, parallelizes model building and minimization
at all levels.

\begin{figure}[htbp]
\begin{center}
\includegraphics[width=7.0cm,height=5cm]{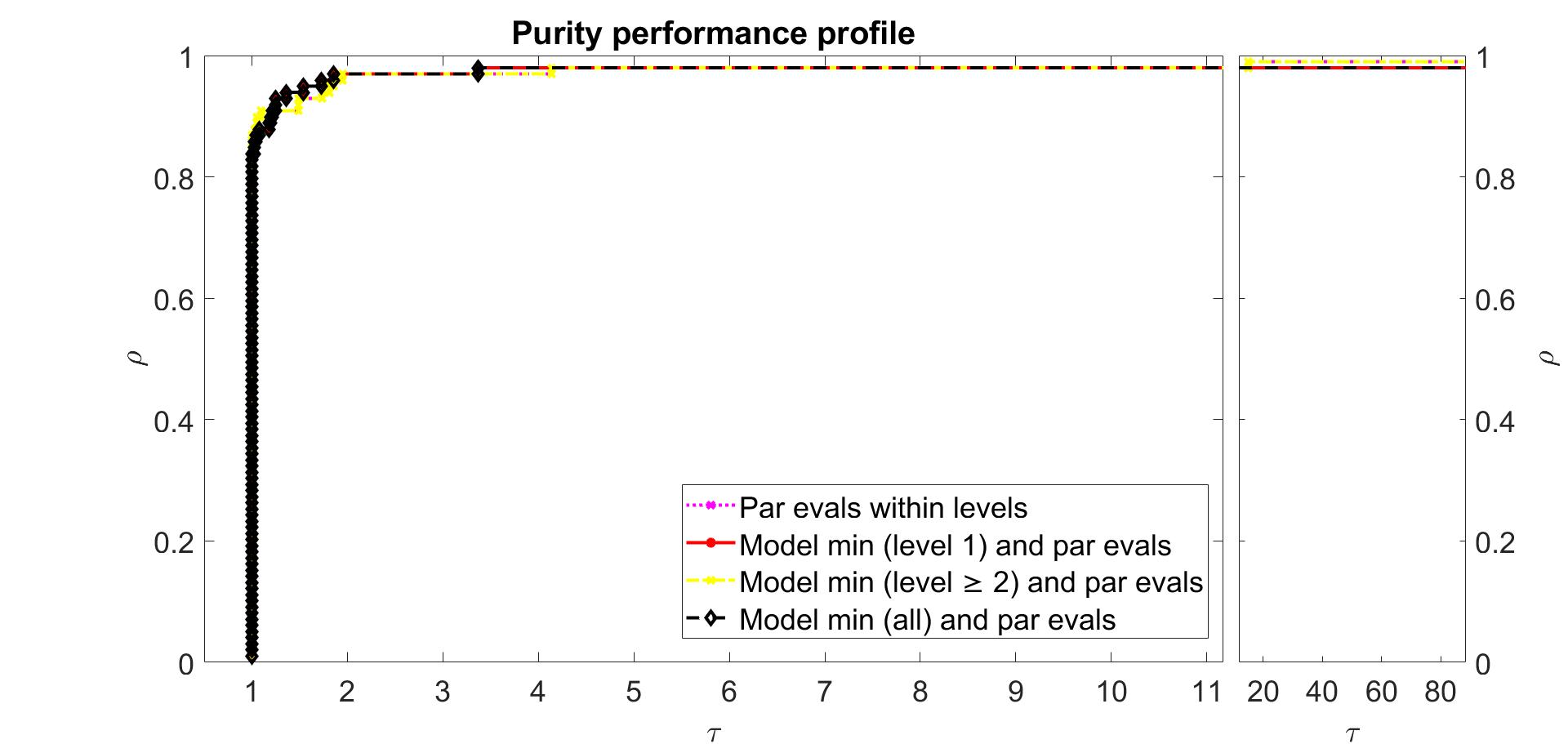}
\includegraphics[width=7.0cm,height=5cm]{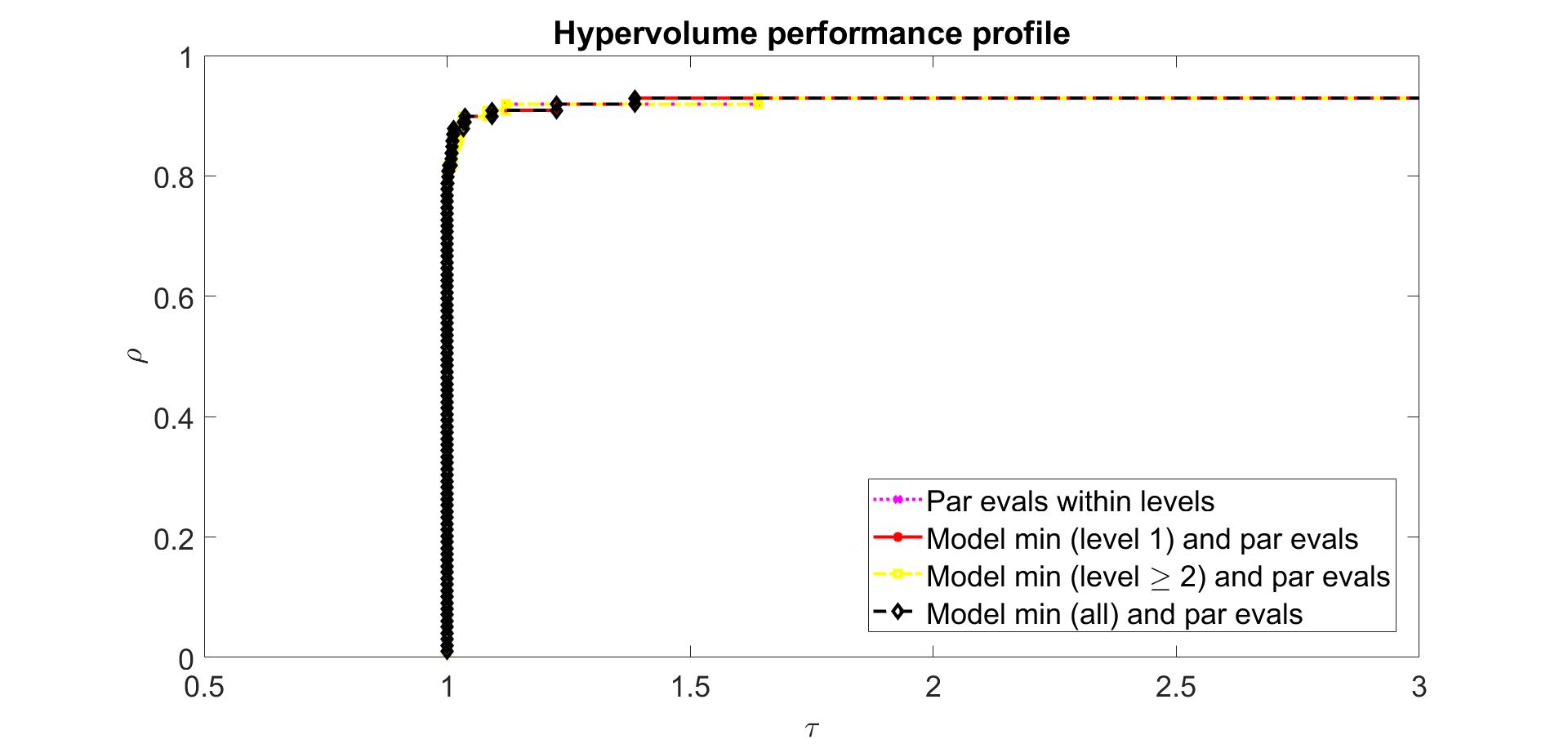}
\end{center}
\caption{\label{fig_min_models_purity_hyper} Performance profiles
corresponding to the parallelization of models building and
minimization, by means of purity and hypervolume metrics.}
\end{figure}

\begin{figure}[htbp]
\begin{center}
\includegraphics[width=7.0cm,height=5cm]{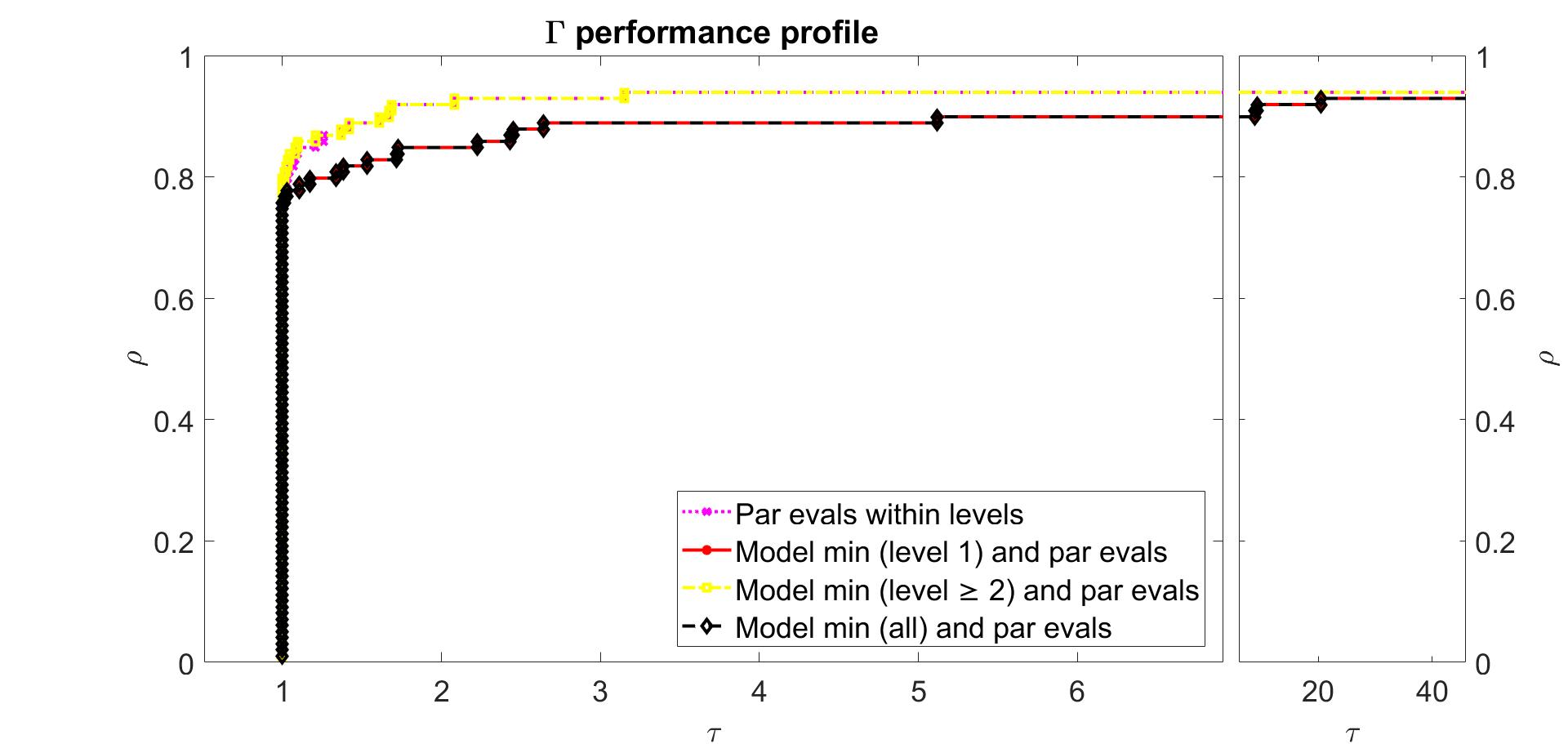}
\includegraphics[width=7.0cm,height=5cm]{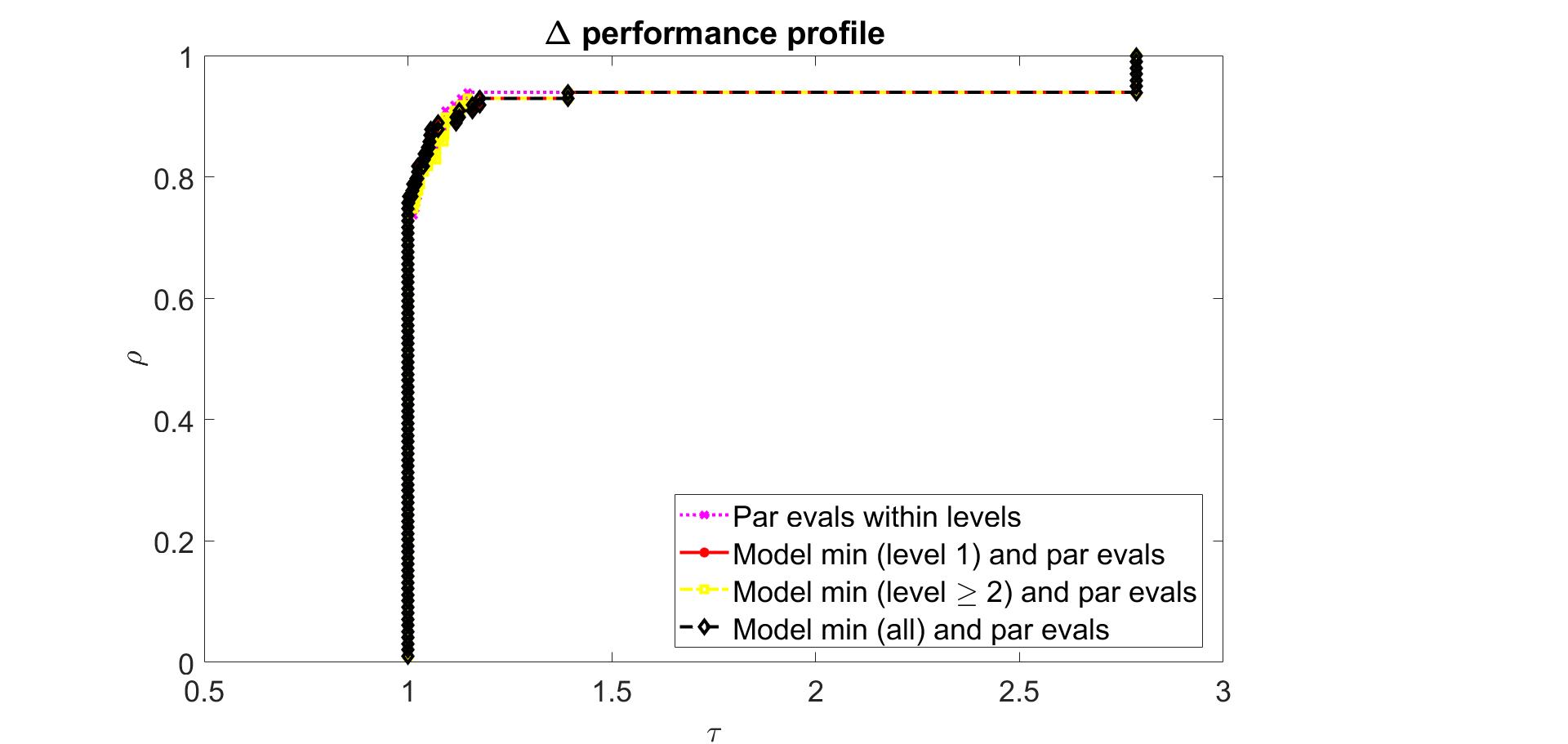}
\end{center}
\caption{\label{fig_min_models_gamma_delta} Performance profiles
corresponding to the parallelization of models building and
minimization, by means of spread metrics ($\Gamma$ and $\Delta$).}
\end{figure}

\begin{figure}[htbp]
\begin{center}
\includegraphics[width=7.0cm,height=5cm]{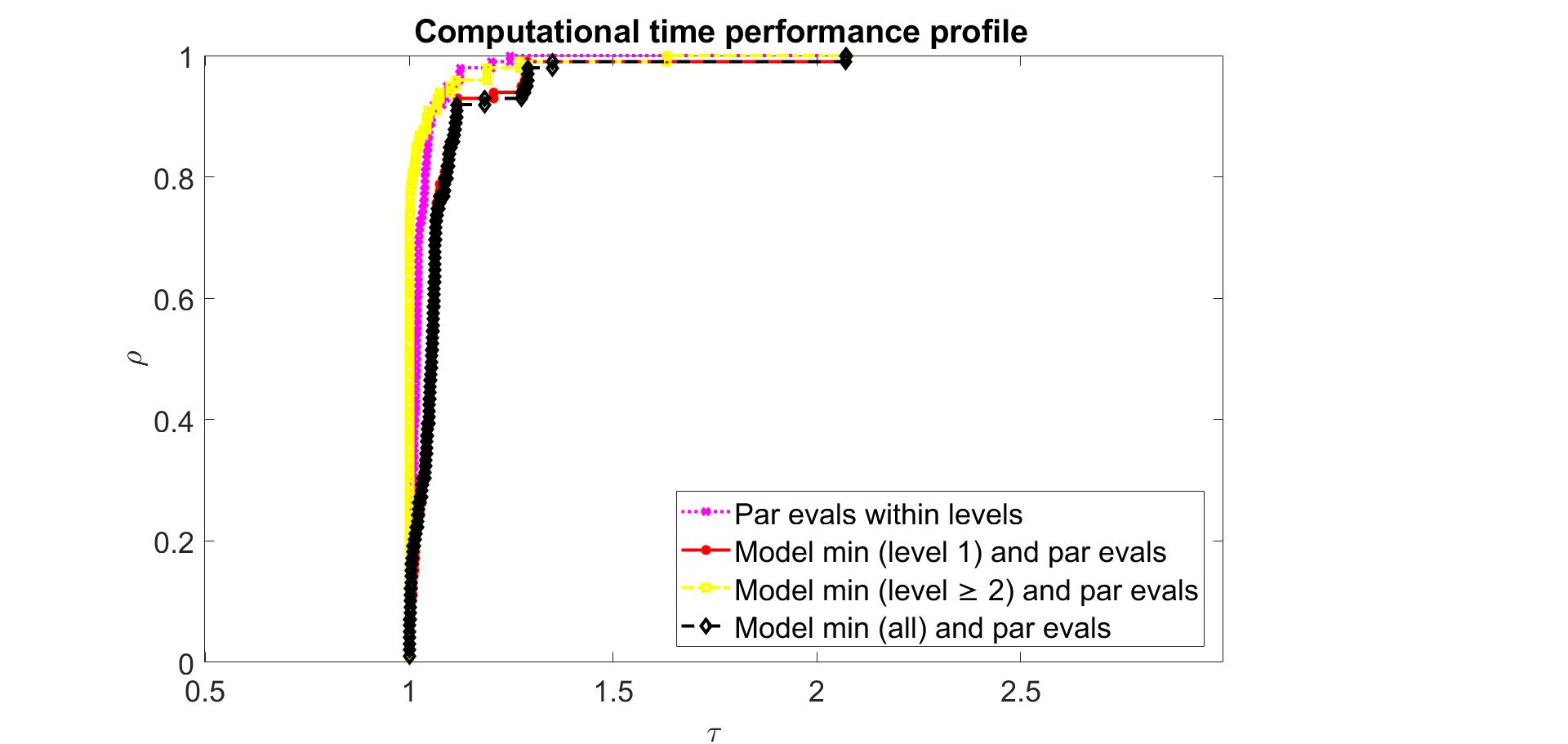}
\end{center}
\caption{\label{fig_min_models_time} Performance profiles
corresponding to the parallelization of models building and
minimization, considering the computational time.}
\end{figure}

Figures~\ref{fig_min_models_purity_hyper},
\ref{fig_min_models_gamma_delta}, and \ref{fig_min_models_time}
report the corresponding performance profiles.

The performance of \emph{Par evals within levels} is very similar
to the one of any of these three new strategies in terms of
purity, hypervolume, and $\Delta$ metrics. The two strategies that
parallelize the model computation present a slightly worse
performance in terms of computational time. This could be
explained by the fact that we are addressing problems with a low
number of components in the objective function, with a reasonably
small dimension ($n\leq30$), which has implications in the linear
systems of equations to be solved for computing the quadratic
models. It is also noteworthy the fact that some of the Matlab's
functions that are used to solve the linear systems of equations
are implicitly parallelized in Matlab. In fact, that justifies the
differences related to the performance of the $\Gamma$ metric,
since the models obtained are different when built sequentially or
in parallel, even when providing exactly the same set of points
for computation.

Since the results did not bring advantages, the strategy \emph{Par
evals within levels} continues to be the default.

\subsection{Iterate point selection based on spread}

As described in Section~\ref{sec:BoostDMS}, each iteration of the
algorithm starts with the selection of an iterate point, which
will be used as model center, at the search step, or, in case of
failure of the search step in finding a new feasible nondominated
point, as poll center, at the poll step. The selection is made
from the list $L_k$ of feasible nondominated points, corresponding
to the point with the largest value of the $\Gamma$
metric~(\ref{gamma_metric}).

This spread metric is considered in an attempt to reduce the gap
between consecutive points lying in the current approximation to
the Pareto front, after projection of each objective function
component in the corresponding dimension. Since the range of each
objective function component can be in very different scales, this
metric could be biased towards one/some of the function
components.

Two additional strategies were considered to prevent this fact.
The first, denoted by \emph{Par evals within levels (Gamma
normalized)}, before computing the largest gap, normalizes the
values of each objective function component $i\in\{1,\ldots,q\}$
using the formula:
$$\displaystyle \frac{f_i(y_j)-\min_{j=1,\ldots,N} f_i(y_j)}{\max_{j=1,\ldots,N} f_i(y_j)-\min_{j=1,\ldots,N}
f_i(y_j)},$$ where $y_1,\ldots,y_N$ represent the points in the
current list of feasible nondominated points $L_k$. This way, a
fairer selection of the largest gap among all components is
expected.

The second strategy, denoted by \emph{Par evals within levels
(Gamma cyclic)}, considers a cyclic approach, changing the
objective function component for which the largest gap is computed
from iteration to iteration in a recurrent way.

\begin{figure}[htbp]
\begin{center}
\includegraphics[width=7.0cm,height=5cm]{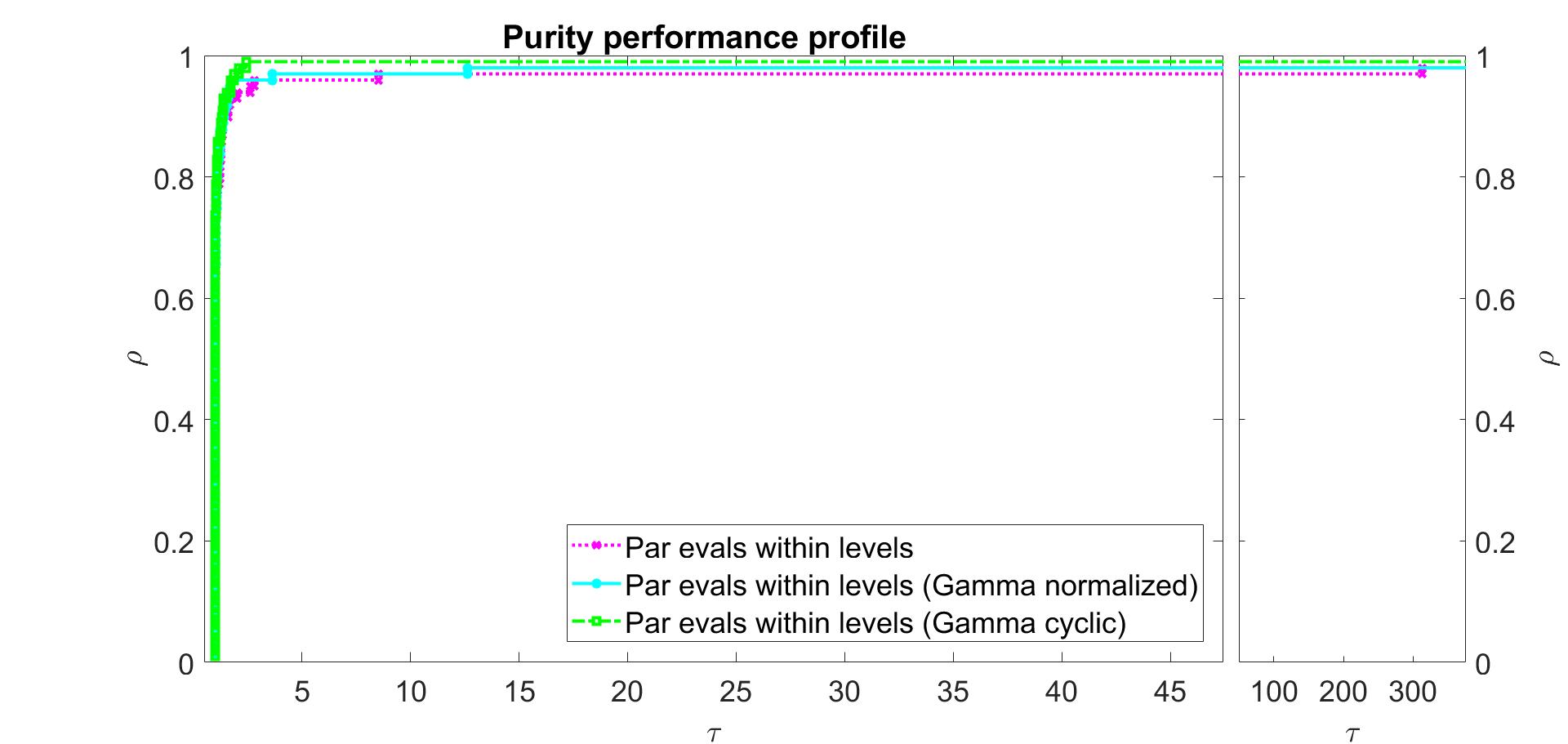}
\includegraphics[width=7.0cm,height=5cm]{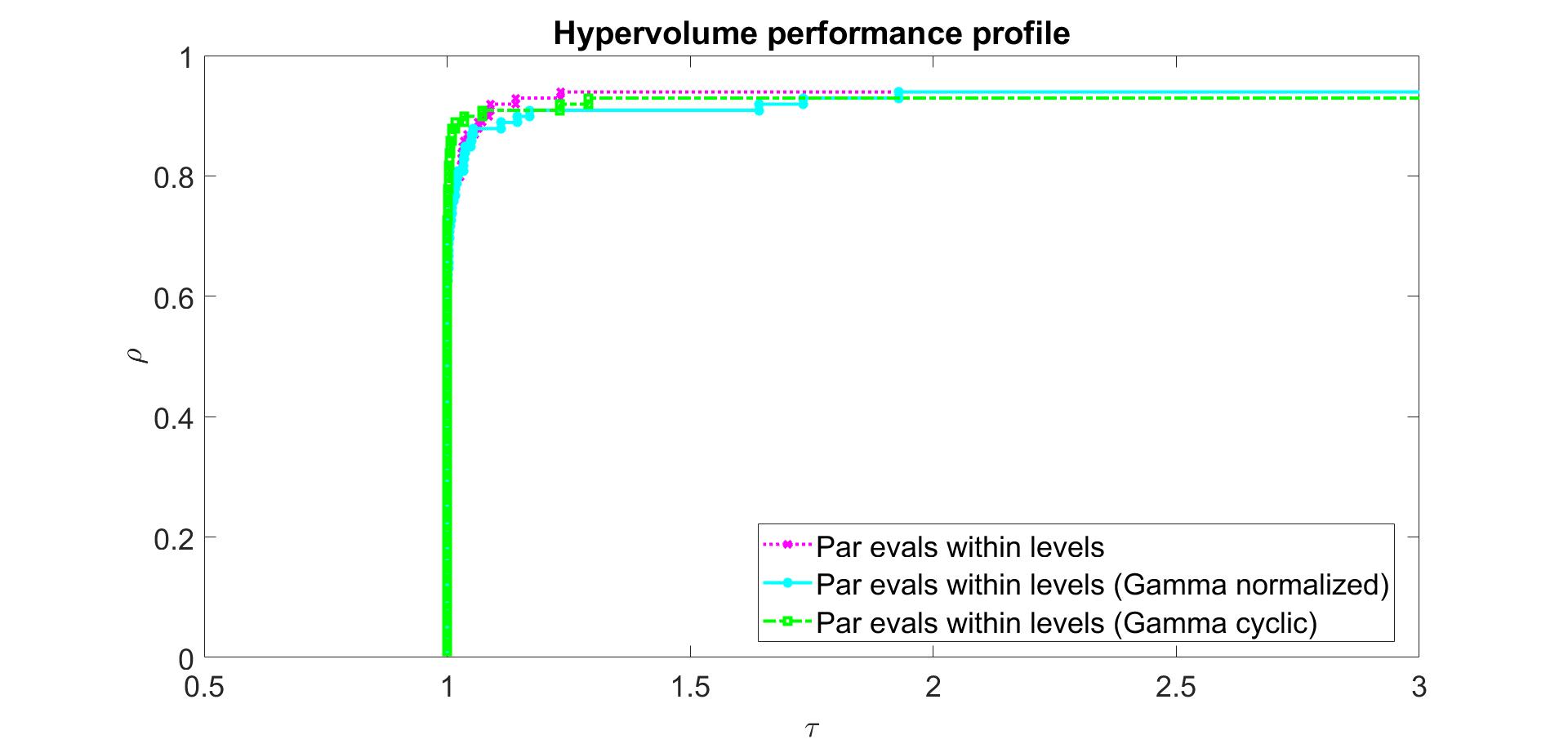}
\end{center}
\caption{\label{fig_spread_purity_hyper} Comparison between
strategies \emph{Par evals within levels}, \emph{Par evals within
levels (Gamma normalized)}, and \emph{Par evals within levels
(Gamma cyclic)}, by means of purity and hypervolume metrics.}
\end{figure}

\begin{figure}[htbp]
\begin{center}
\includegraphics[width=7.0cm,height=5cm]{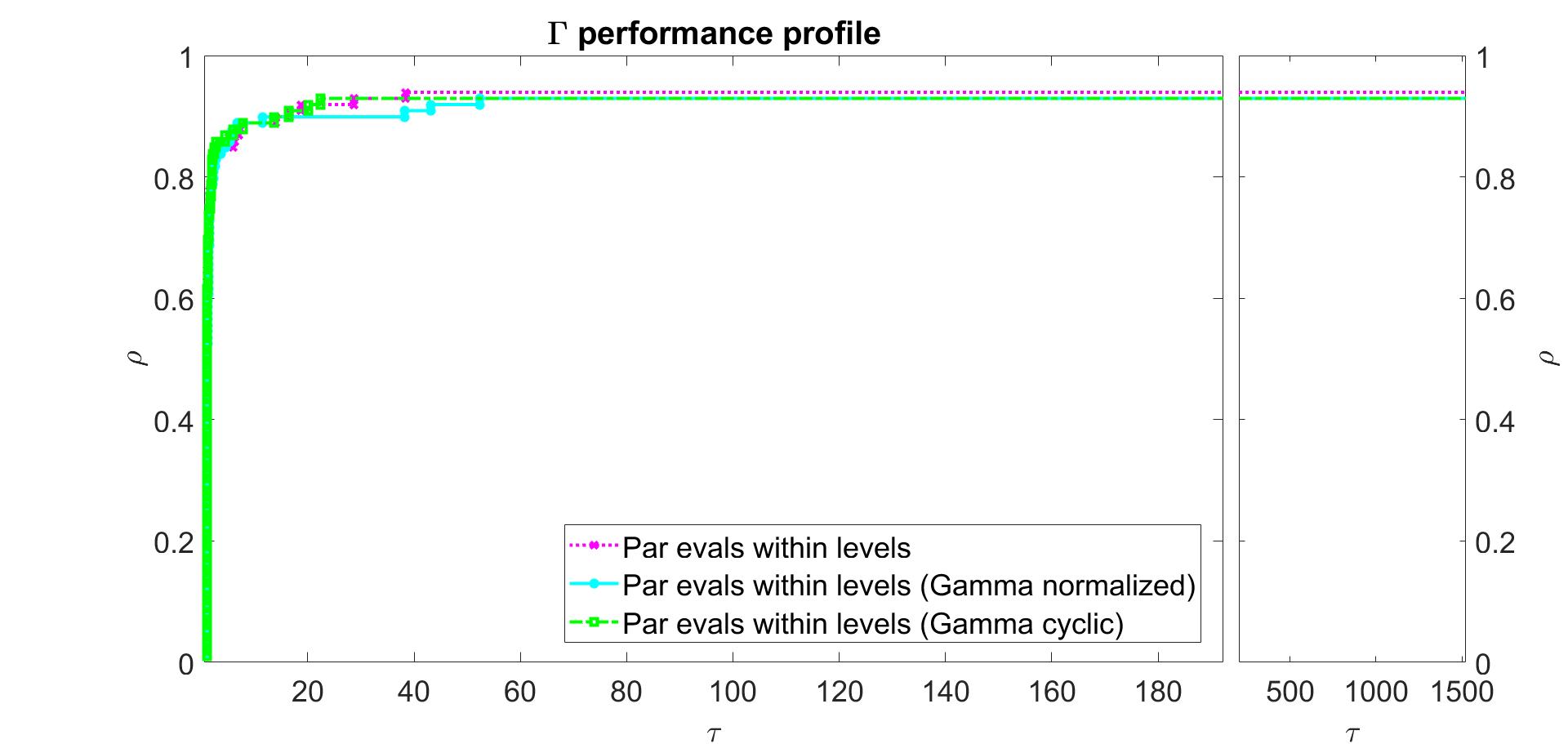}
\includegraphics[width=7.0cm,height=5cm]{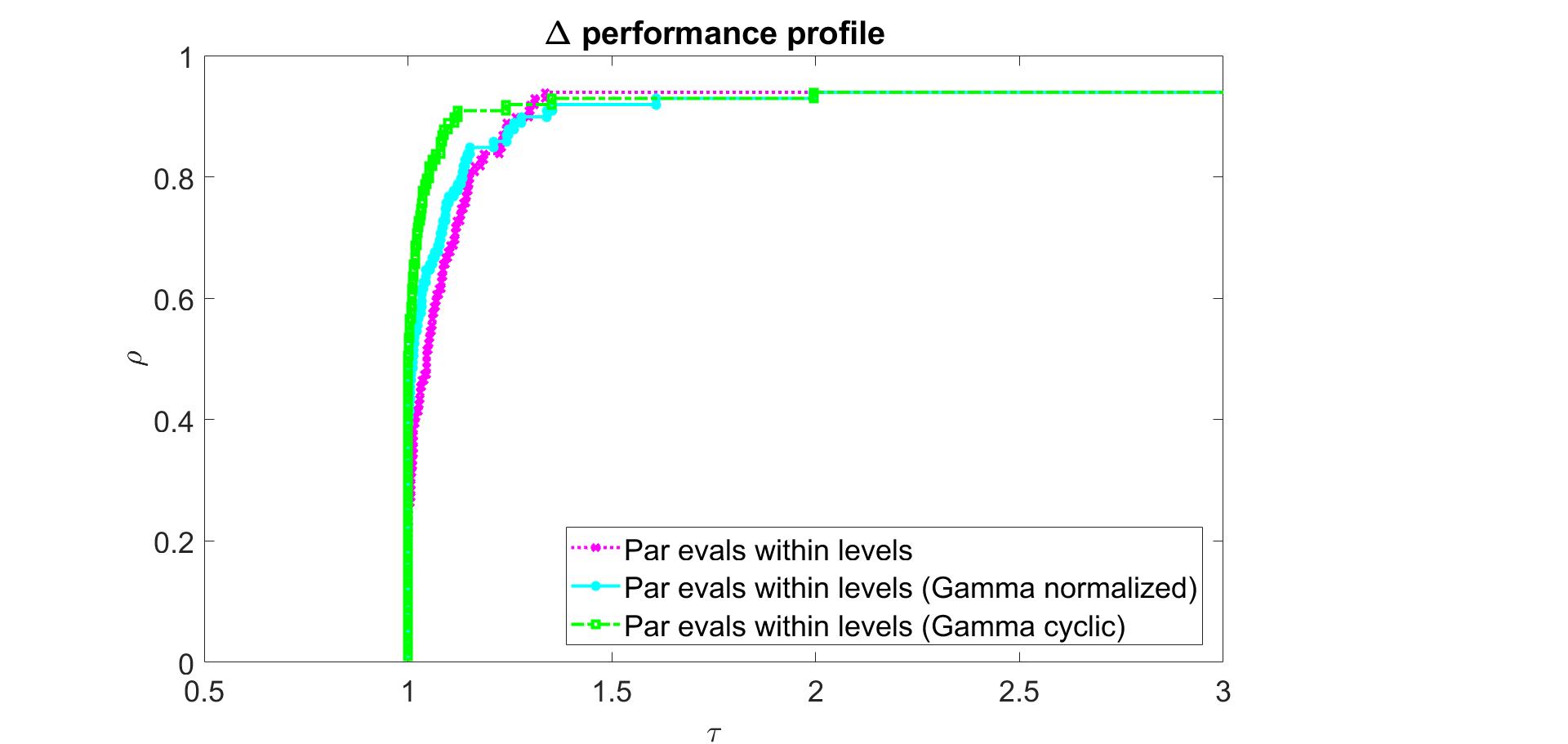}
\end{center}
\caption{\label{fig_spread_gamma_delta} Comparison between
strategies \emph{Par evals within levels}, \emph{Par evals within
levels (Gamma normalized)}, and \emph{Par evals within levels
(Gamma cyclic)}, by means of spread metrics ($\Gamma$ and
$\Delta$).}
\end{figure}

\begin{figure}[htbp]
\begin{center}
\includegraphics[width=7.0cm,height=5cm]{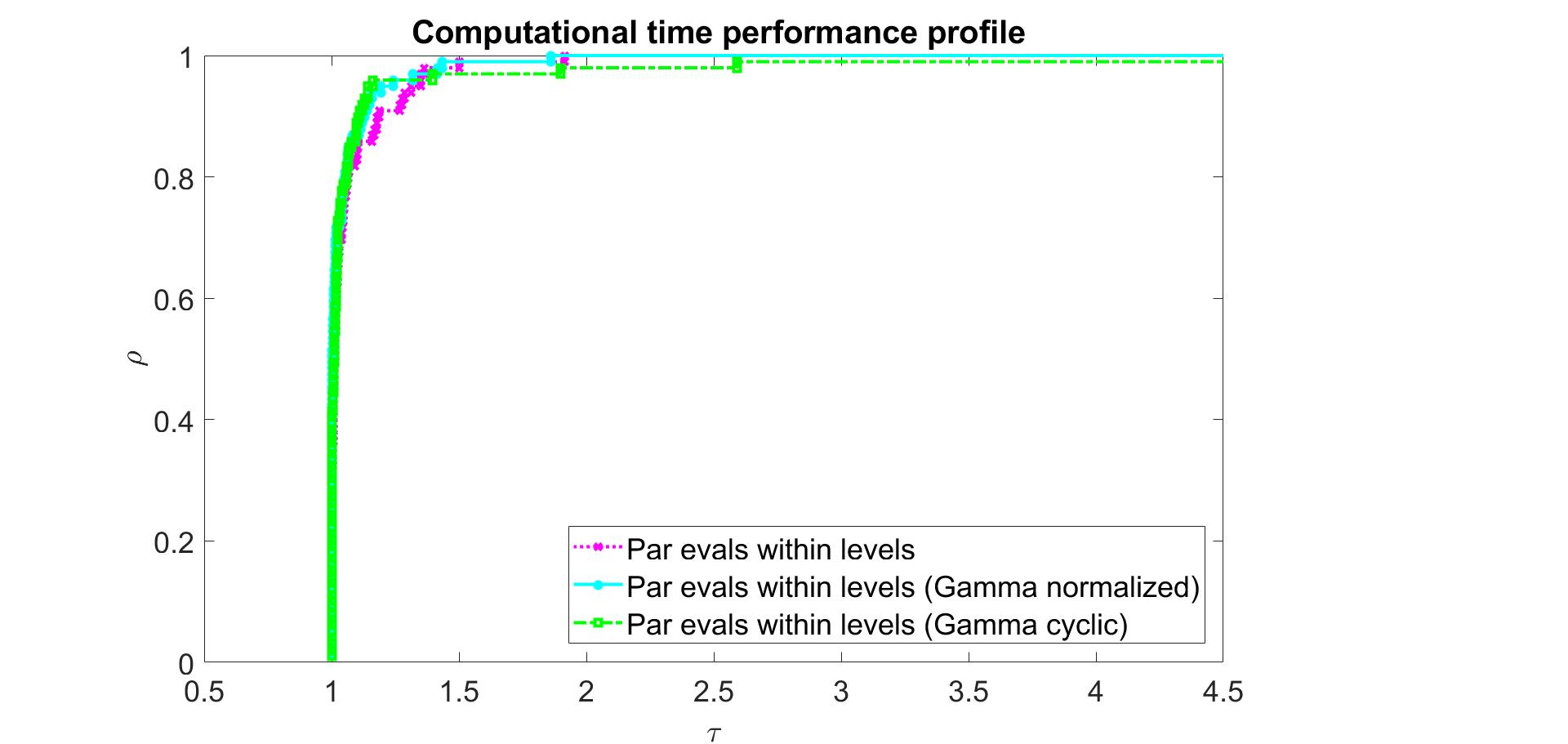}
\end{center}
\caption{\label{fig_spread_time}Comparison between strategies
\emph{Par evals within levels}, \emph{Par evals within levels
(Gamma normalized)}, and \emph{Par evals within levels (Gamma
cyclic)}, considering the computational time.}
\end{figure}

Figures~\ref{fig_spread_purity_hyper},
\ref{fig_spread_gamma_delta}, and \ref{fig_spread_time}  compare
these two strategies with \emph{Par evals within levels}. As
usual, function evaluation is parallelized both at the search and
the poll steps, in the former respecting the levels structure.

Although the results are very similar for hypervolume and $\Gamma$
metrics, a slight advantage of \emph{Par evals within levels
(Gamma cyclic)} is noticed for purity, $\Delta$ and computational
time, which justifies our option to use it as strategy for the
iterate point selection.

\subsection{Selection of more than one iterate point}
Iterate points are important not only as base points for trying to
close gaps in the Pareto front, but can also be a way of expanding
the Pareto front, in an attempt of reaching the corresponding
extreme points. Several parallel strategies were developed,
focusing on these two goals. Variants included the selection of
more than one iterate point per iteration, considering the
$\Gamma$ metric applied in a cyclic way, the objective function
value (minimum and maximum values at each component), and the
stepsize parameter (largest indicating more promising points). We
only report results for the two most successful variants.

Both new strategies take into account the two points defining the
largest gap for the objective function component corresponding to
the current iteration as well as the $q$ points corresponding to
the lowest value for each objective function component. In the
first strategy, denoted by \emph{It centers based on spread $+$
best values}, these $q + 2$ points are used as iterate points. In
the second strategy, \emph{It centers based on spread $\vee$ best
values}, iterations alternate between the selection of the two
points corresponding to the largest gap and the $q$ points
corresponding to the lowest objective function components values.
As usual, function evaluation is performed in parallel, keeping
the levels structure at the search step.

\begin{figure}[htbp]
\begin{center}
\includegraphics[width=7.0cm,height=5cm]{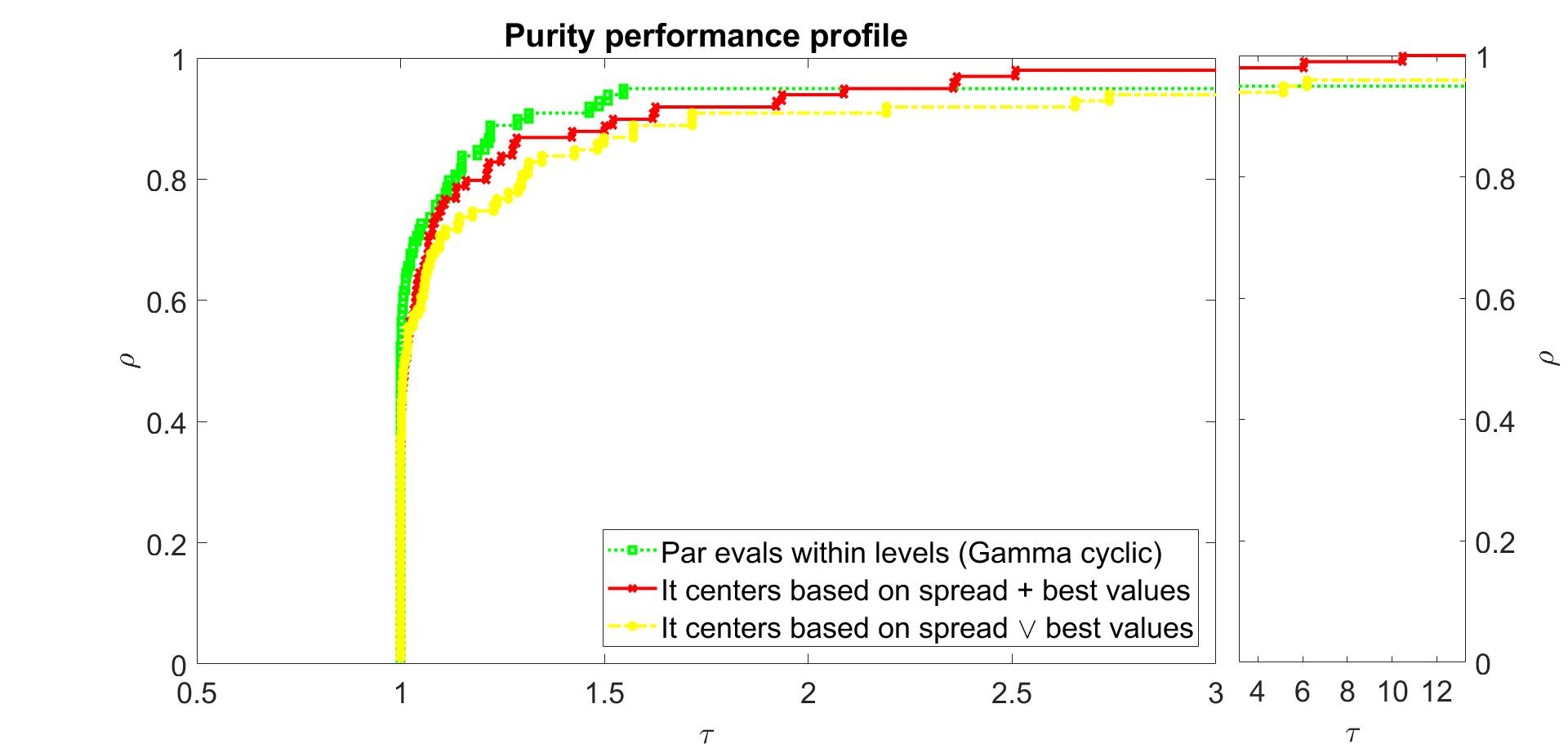}
\includegraphics[width=7.0cm,height=5cm]{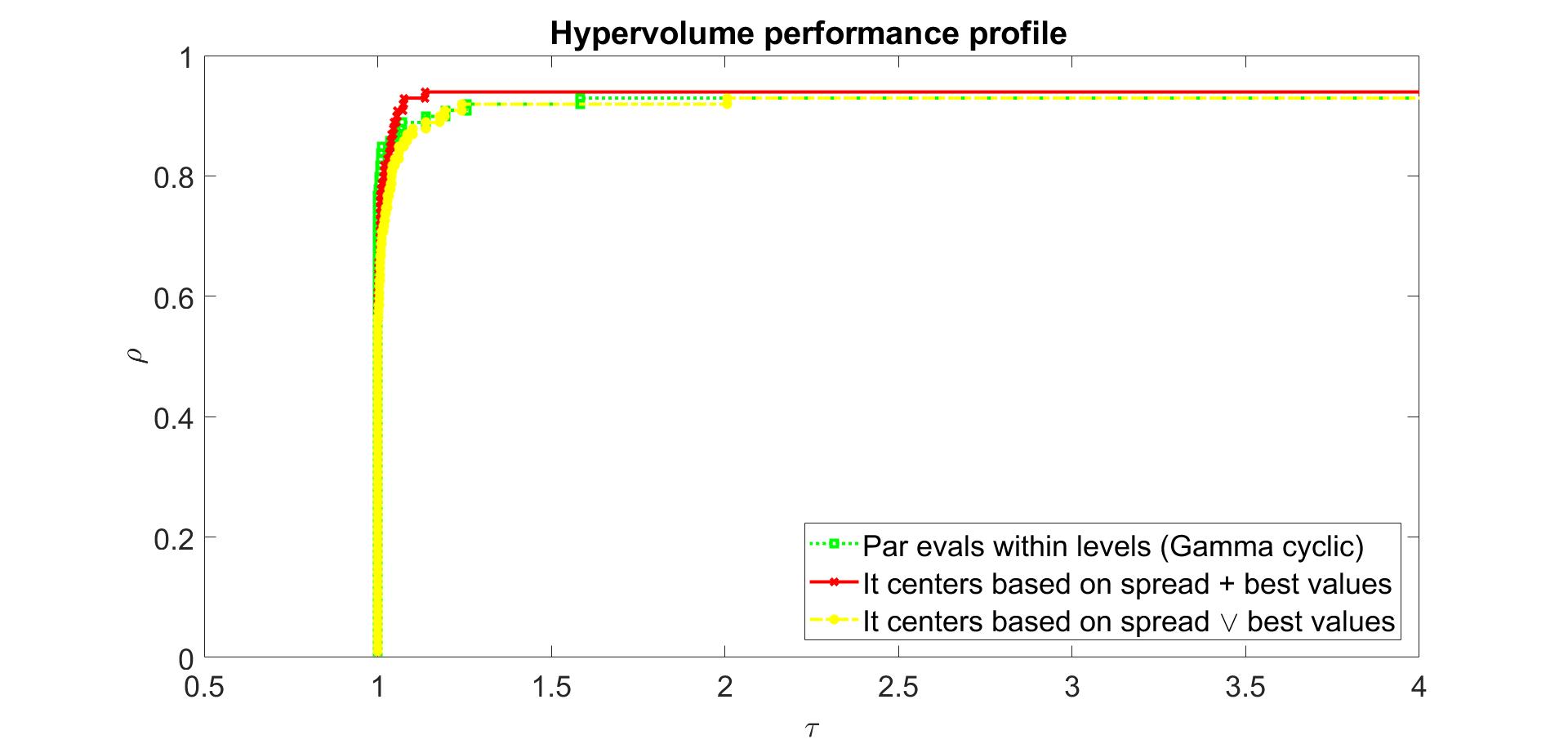}
\end{center}
\caption{\label{fig_best_values_purity_hyper} Comparison between
parallel strategies \emph{Par evals within levels (Gamma cyclic)},
\emph{It centers based on spread $+$ best values}, and \emph{It
centers based on spread $\vee$ best values}, by means of purity
and hypervolume metrics.}
\end{figure}

\begin{figure}[htbp]
\begin{center}
\includegraphics[width=7.0cm,height=5cm]{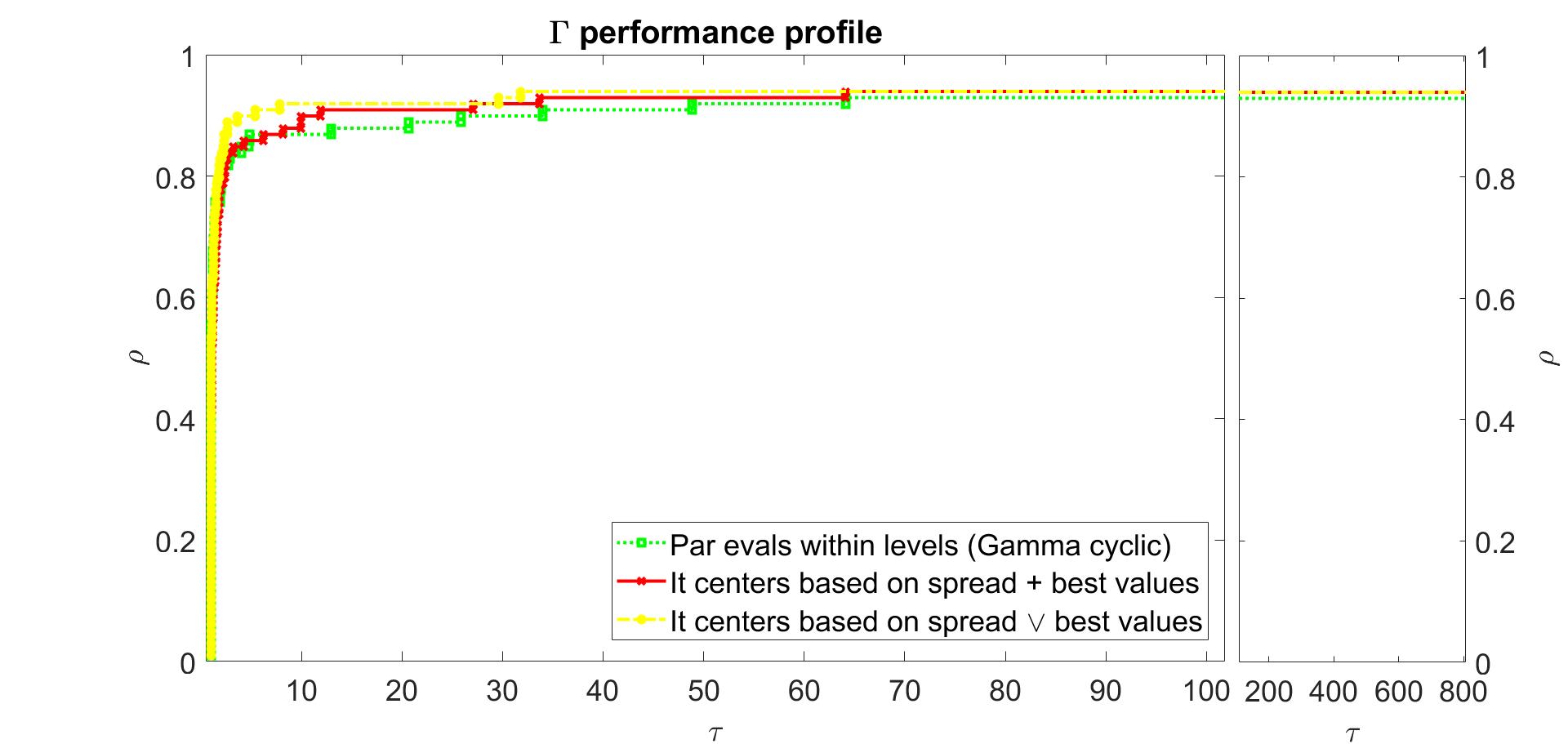}
\includegraphics[width=7.0cm,height=5cm]{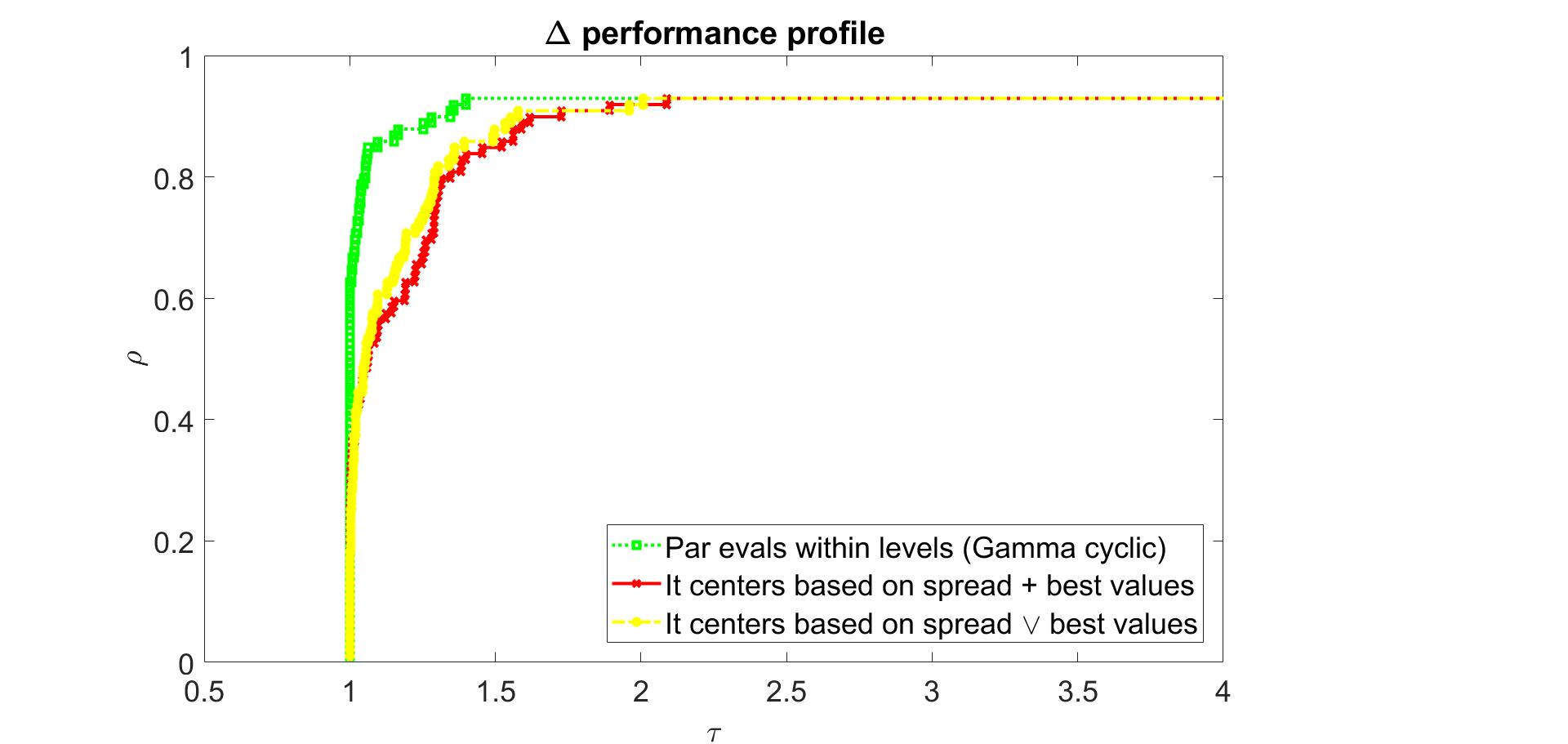}
\end{center}
\caption{\label{fig_best_values_gamma_delta} Comparison between
parallel strategies \emph{Par evals within levels (Gamma cyclic)},
\emph{It centers based on spread $+$ best values}, and \emph{It
centers based on spread $\vee$ best values}, by means of spread
metrics ($\Gamma$ and $\Delta$).}
\end{figure}

\begin{figure}[htbp]
\begin{center}
\includegraphics[width=7.0cm,height=5cm]{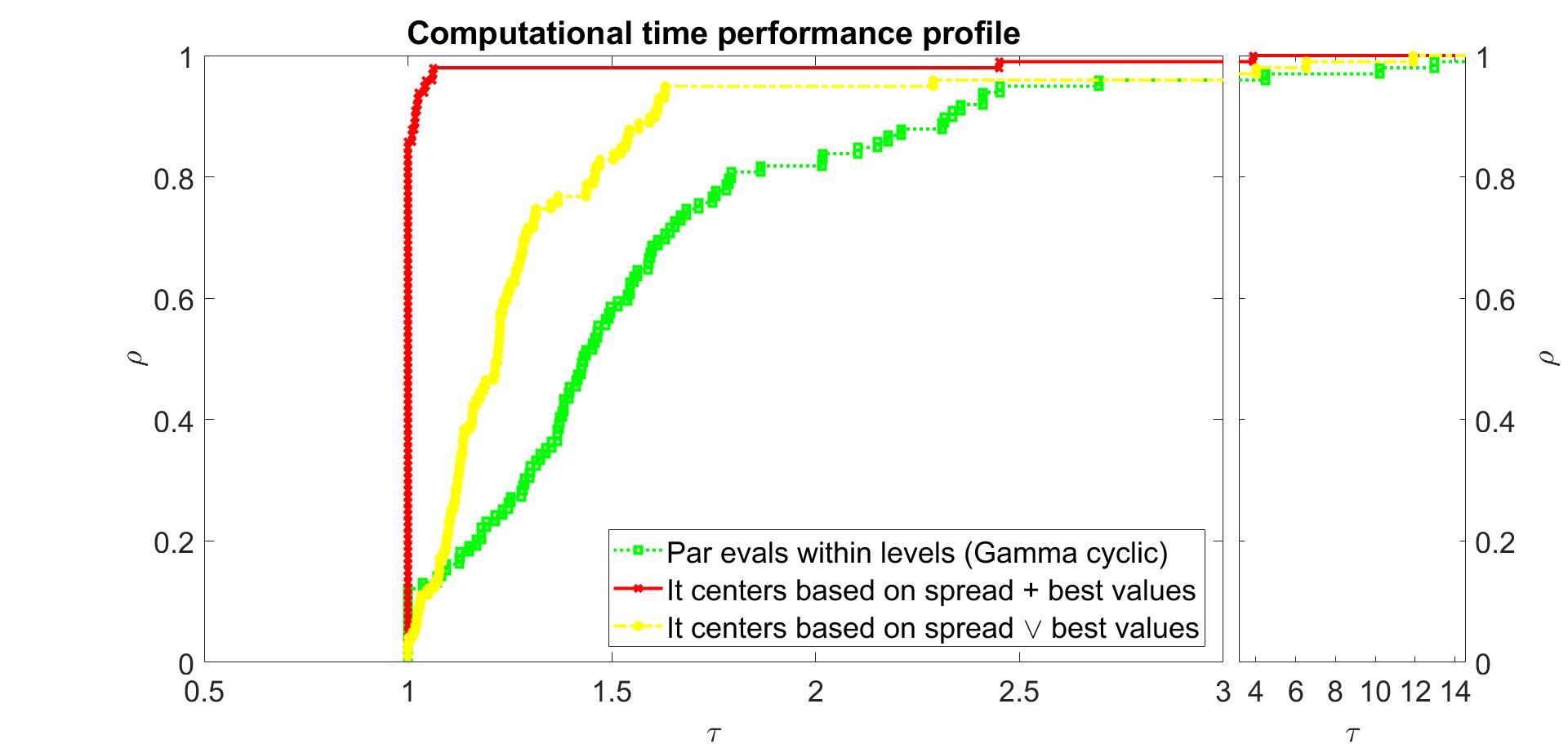}
\end{center}
\caption{\label{fig_best_values_time} Comparison between parallel
strategies \emph{Par evals within levels (Gamma cyclic)}, \emph{It
centers based on spread $+$ best values}, and \emph{It centers
based on spread $\vee$ best values}, considering the computational
time.}
\end{figure}

Figures~\ref{fig_best_values_purity_hyper},
\ref{fig_best_values_gamma_delta}, and \ref{fig_best_values_time}
compare these two strategies with \emph{Par evals within levels
(Gamma cyclic)}.

Although \emph{Par evals within levels (Gamma cyclic)} presents a
better performance for purity and $\Delta$ metrics, strategy
\emph{It centers based on spread $+$ best values} clearly
outperforms the other strategies in computational time, and
presents a slightly better performance in what respects to the
hypervolume metric, justifying our option for it.

\subsection{Selected parallelization approach}\label{sec:best_BoostDMS}

Finally, we compared our selected parallel strategy with the
sequential variant of BoostDMS. At each iteration, this strategy
considers $q+2$ iterate points from the current list $L_k$,
corresponding to the $q$ points with the lowest value for each
component of the objective function and to the two points defining
the largest gap in the component of the objective function
associated to the current iteration (which changes in a cyclic way
between iterations). These $q+2$ points will be used as model
centers at the search step and all those that fail in generating a
new feasible nondominated point will be used as poll centers.
Function evaluations are always performed in parallel, at the
search and poll steps. In the former, the levels structure is
respected, only evaluating a given level if the previous one was
unable to generate a new feasible nondominated point.

\begin{figure}[htbp]
\begin{center}
\includegraphics[width=7.0cm,height=5cm]{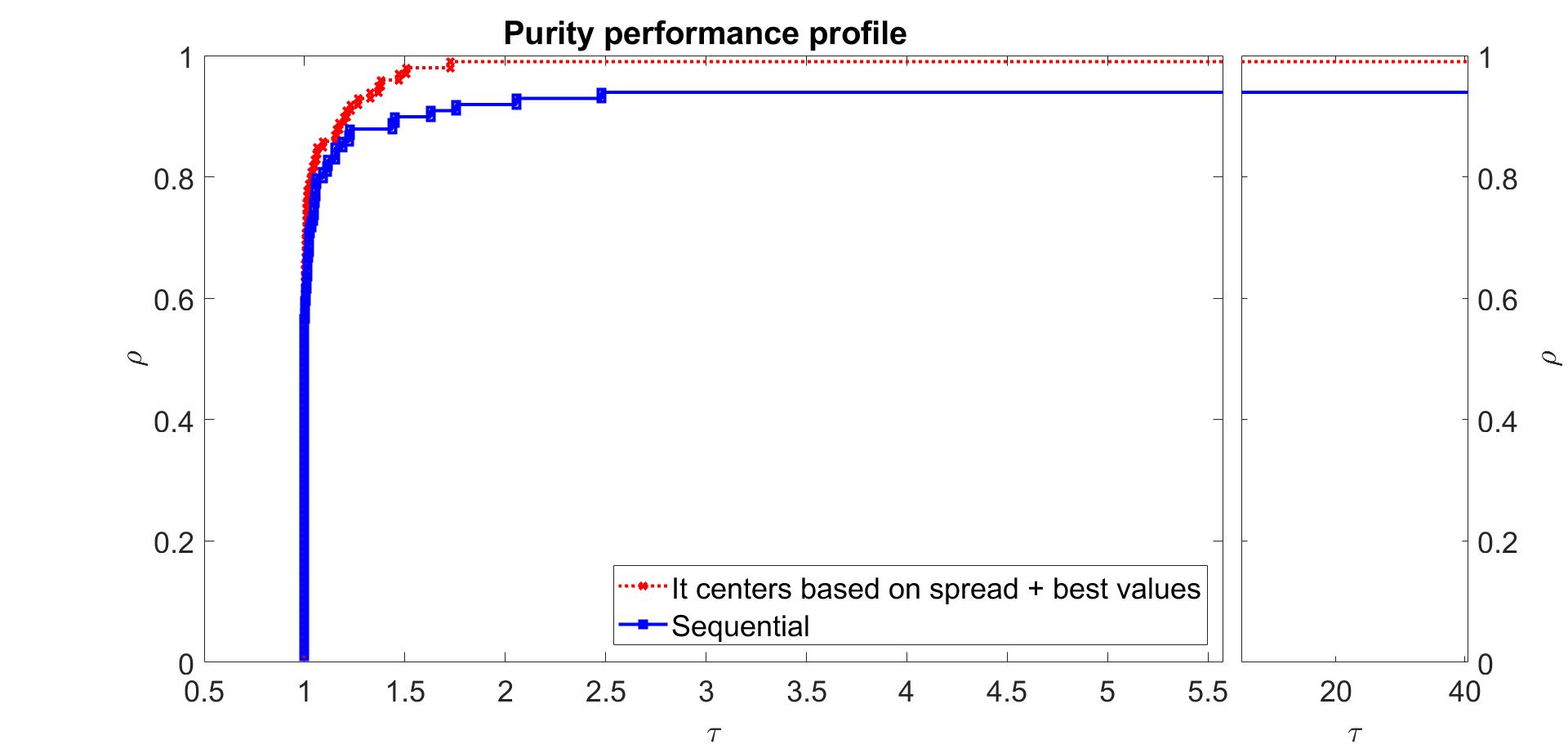}
\includegraphics[width=7.0cm,height=5cm]{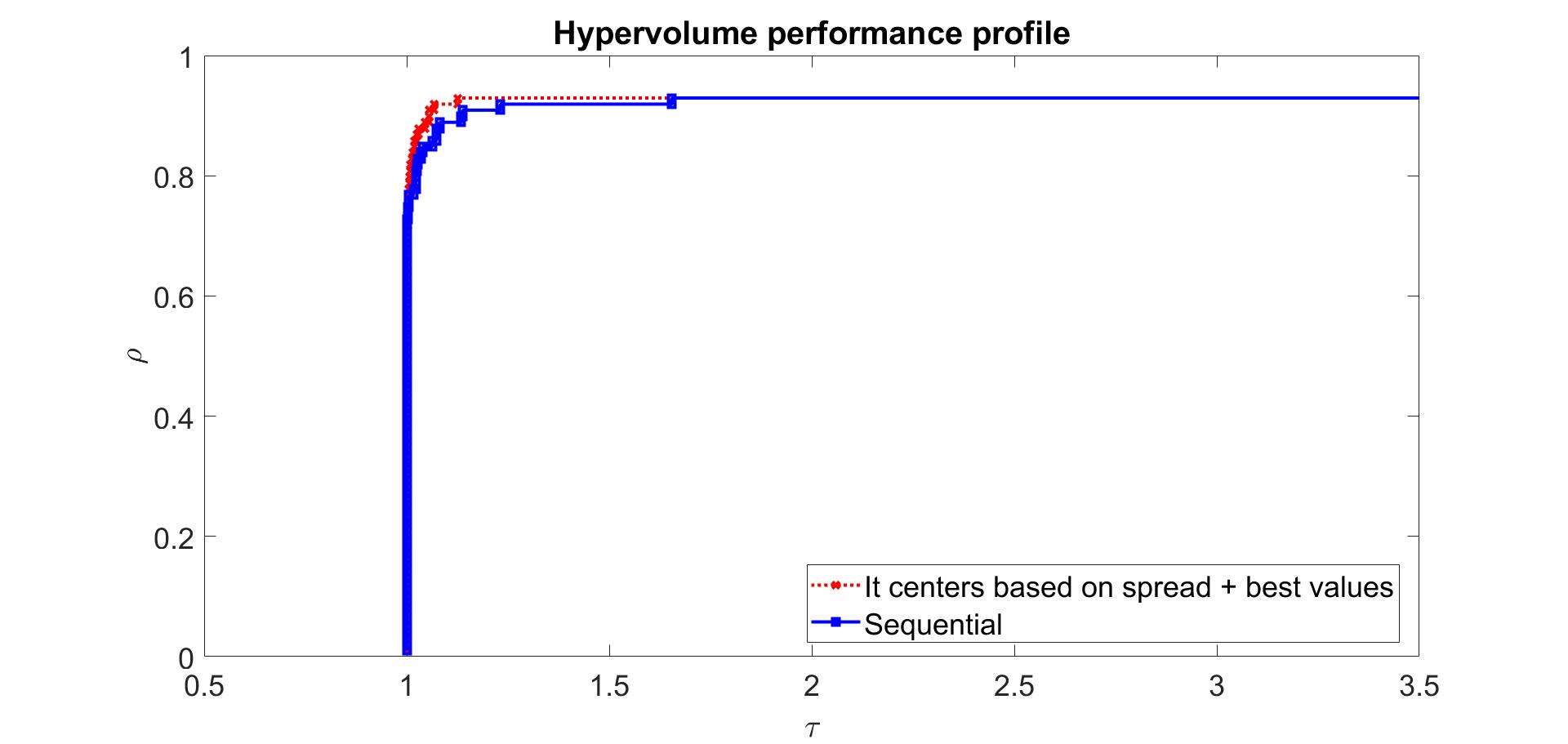}
\end{center}
\caption{\label{fig_final_purity_hyper} Comparison between the
selected parallel version and the sequential implementation of
BoostDMS, by means of purity and hypervolume metrics.}
\end{figure}

\begin{figure}[htbp]
\begin{center}
\includegraphics[width=7.0cm,height=5cm]{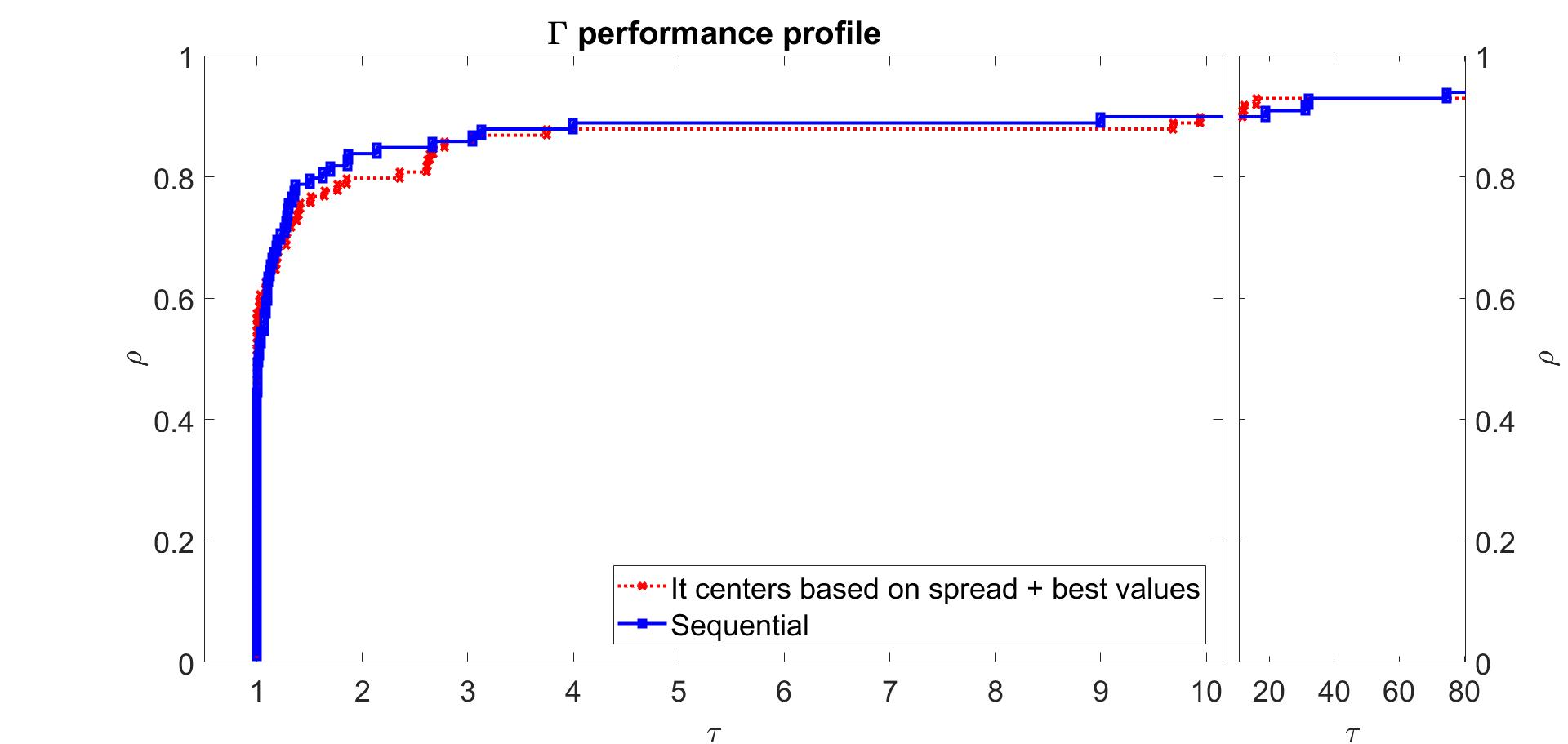}
\includegraphics[width=7.0cm,height=5cm]{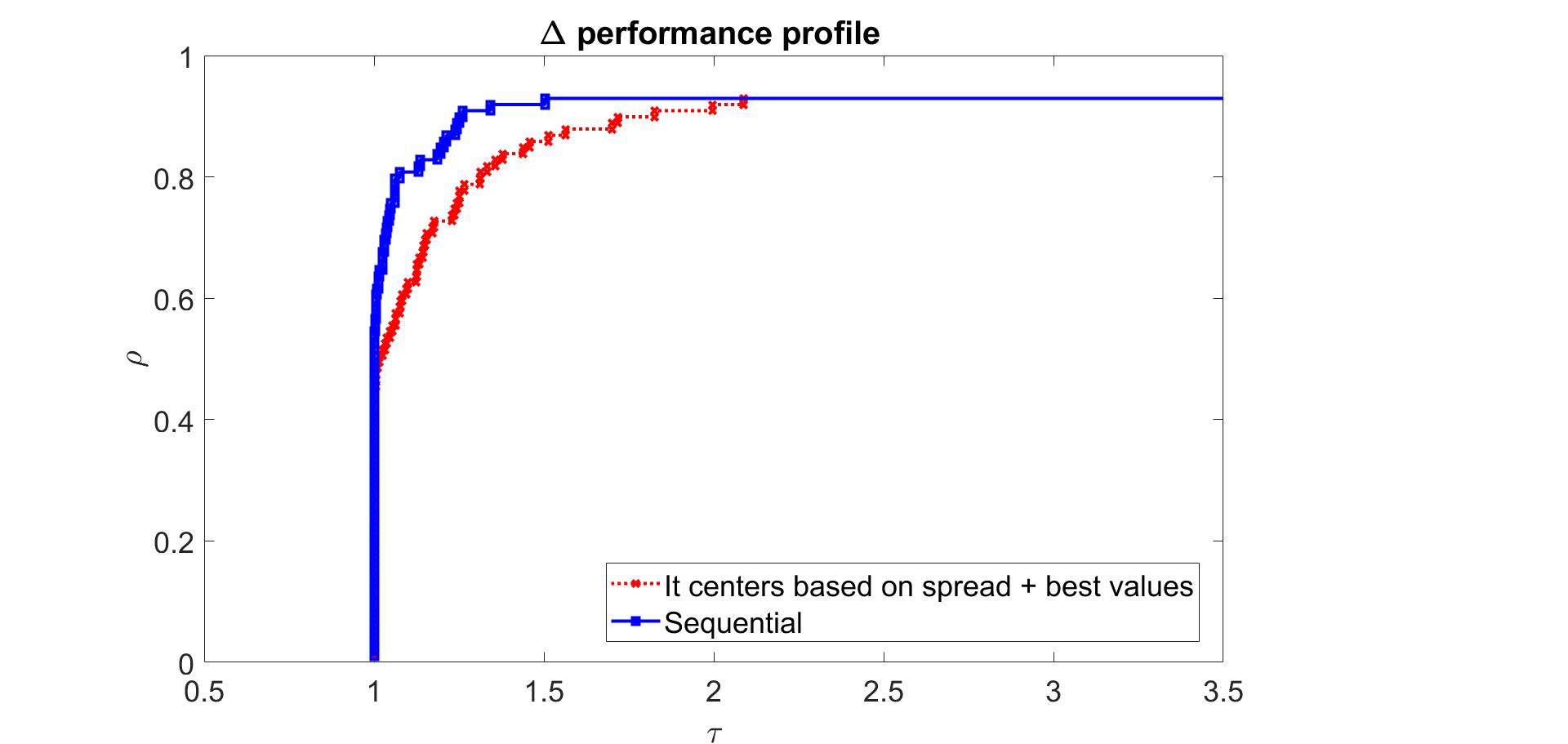}
\end{center}
\caption{\label{fig_final_gamma_delta} Comparison between the
selected parallel version and the sequential implementation of
BoostDMS, by means of spread metrics ($\Gamma$ and $\Delta$).}
\end{figure}

\begin{figure}[htbp]
\begin{center}
\includegraphics[width=7.0cm,height=5cm]{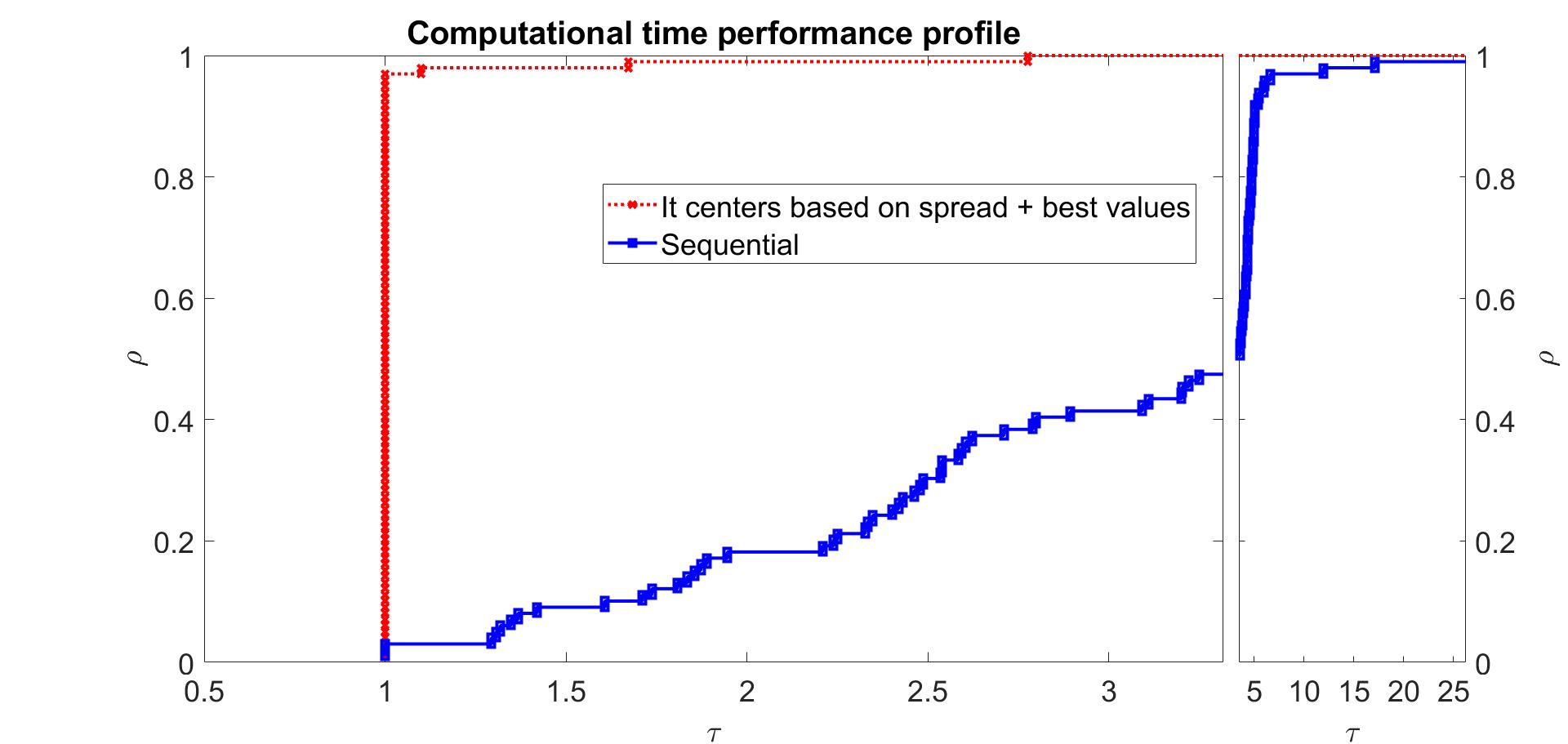}
\end{center}
\caption{\label{fig_final_time} Comparison between the selected
parallel version and the sequential implementation of BoostDMS,
considering computational time.}
\end{figure}

Figures~\ref{fig_final_purity_hyper}, \ref{fig_final_gamma_delta},
and \ref{fig_final_time} report the performance profiles comparing
\emph{It centers based on spread $+$ best values} strategy with
the initial sequential version of BoostDMS.

Although there is a decrease in performance for spread (see
Figure~\ref{fig_final_gamma_delta}), the advantage brought by the
parallel strategy selected is clear. Not only it presents a huge
increase in performance regarding computational time
(Figure~\ref{fig_final_time}), but also a better performance for
purity and a slight advantage in terms of hypervolume (see
Figure~\ref{fig_final_purity_hyper}).

\section{A chemical engineering application}\label{sec:styrene}

Often academic problems, like the ones considered in
Section~\ref{sec:numerical_results}, do not reflect all the
challenges of real applications. In this section we consider the
real chemical engineering problem, related to the production of
styrene, as described in~\cite{Audet_et_al_2008,Audet_et_al_2010}.

The styrene production process involves four steps: preparation of
reactants, catalytic reactions, a first distillation -- where
styrene is recovered, and a second distillation -- where benzene
is recovered. During this second distillation, unreacted
ethylbenzene is recycled, as an initial reactant on the process.
This is a tri-objective problem ($q=3$), where it is intended to
maximize the net present value associated to the process
($f_{1}$), the purity of the produced styrene ($f_{2}$) and the
overall ethylbenzene conversion into styrene
($f_{3}$)~\cite{Audet_et_al_2010}. The problem has 8 variables,
subject to bounds constraints, and 9 other constraints, some
process-related (\textit{e.g.} environmental norms regarding
excesses) and some economical (\textit{e.g.} investment value).
More details can be found in~\cite{Audet_et_al_2008}.

A C++ numerical simulator has been developed for this chemical
engineering process, where each complete simulation corresponds to
a function evaluation. Due to the presence of recycling loops,
like the one of ethylbenzene, the process must be entirely
simulated until a result is provided. The total time associated to
the computation of a function value fluctuates, with an average of
1 second, when the code succeeds in evaluating a given point.
Additionally, the problem reveals hidden constraints, since
feasible points, according to the defined constraints, often fail
to produce a finite numerical value for the objective function. In
this case, the simulation is generally faster, with an average
computational time of 0.001 seconds. Note that these computational
times were obtained with the hardware/software configuration
described in Section~\ref{sec:numerical_results}.

\begin{figure}[htbp]
\begin{center}
\includegraphics[width=10cm,height=6cm]{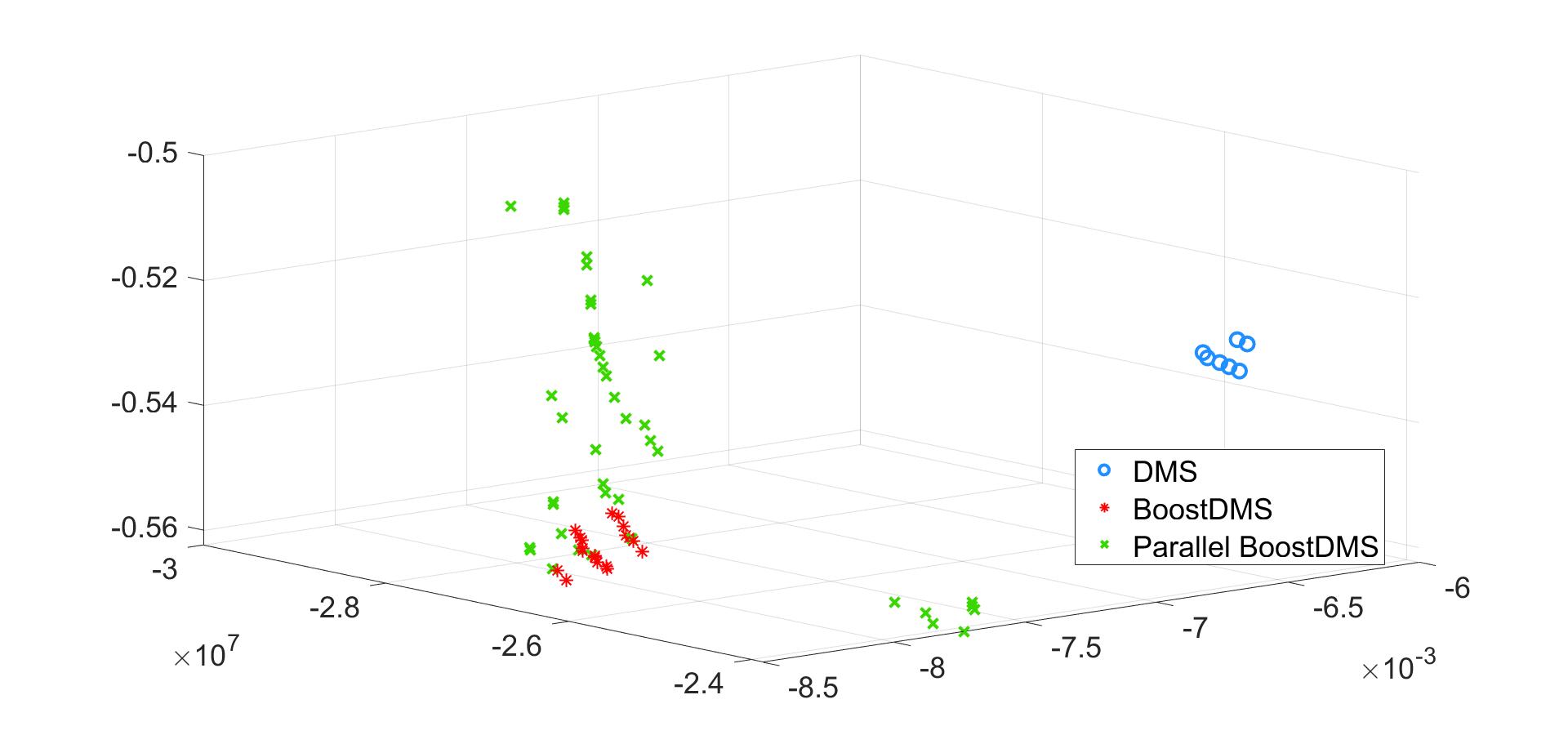}
\end{center}
\caption{\label{fig_styrene} Final approximations to the Pareto
front computed for the chemical engineering problem by solvers
DMS, BoostDMS, and the selected parallel version of BoostDMS.}
\end{figure}

We used DMS~\cite{ALCustodio_et_al_2011},
BoostDMS~\cite{CBras_ALCustodio_2020} and the selected parallel
version of BoostDMS, described in Section~\ref{sec:best_BoostDMS},
to solve the problem. All the defaults were considered for the
solvers, with exception of the initialization, where line sampling
was replaced by a single point, provided in
literature~\cite{Audet_et_al_2008}. Figure~\ref{fig_styrene}
corresponds to the approximations to the Pareto front generated by
each one of the three solvers.

The 7 points in the final solution computed by DMS are all
dominated by each one of the approximations to the Pareto front
computed by BoostDMS or its parallel version (which comprise 19
and 44 nondominated points, respectively). It is also clear the
gain in volume of the dominated region obtained with both versions
of BoostDMS. Table~\ref{tab_styrene} presents the values of the
different multiobjective metrics for the three solvers, as well as
the computational time (an average of five runs) and the total
number of function evaluations performed.

\begin{table}
\begin{scriptsize}
\begin{tabular}{c c c c c c c}
\hline \textbf{Solver} & \textbf{Purity} & \textbf{Hypervolume} &
\textbf{Gamma}  & \textbf{Delta}   & \textbf{AvgTime(s)} &
\textbf{FuncEvals} \\ \hline
DMS & 0\% & 0.002 & 3.96e+06 & 1.071 & 3987 & 9260 \\
BoostDMS & 53\% & 0.70 & 3.13e+06 & 1.173 & 6970 & 14978 \\
Parallel BoostDMS & 89\% & 0.78 & 2.52e+06 & 1.224 & 1670 &
20000\\\hline
\end{tabular}
\end{scriptsize}
\caption{\label{tab_styrene} Metrics associated to the solution of
the chemical engineering problem, computed by solvers DMS,
BoostDMS, and the selected parallel version of BoostDMS.}
\end{table}

BoostDMS is considerably slower than DMS, but presents a
remarkable improvement in the quality of the solution computed.
This quality increases in the solution computed by the parallel
version of BoostDMS, particularly in what respects to purity.
These results corroborate what was observed in the academic test
set, with the parallel version of BoostDMS allowing relevant gains
in terms of purity, a slight advantage regarding hypervolume, but
a worse performance in spread metric $\Delta$.

Although the parallel version of BoostDMS performs a larger number
of function evaluations than any of the two other solvers, it is
considerably faster than any of the two sequential versions. In
fact, if we compute the speedup associated to the parallel
version, we will get a value of $4.17$, corresponding to an
efficiency of $52.17\%$ (for the test set of
Section~\ref{sec:numerical_results}, the average speedup is in the
interval $[3.22;4.45]$ and the average efficiency in
$[40.31\%;55.68\%]$, in both cases with $95\%$ of confidence).
However, considering that the number of function evaluations
performed by the sequential version is clearly lower than the one
of the parallel version, a crude correction, assuming that a total
of $20\,000$ function evaluations was performed by both variants,
would give a speedup of $5.6$, corresponding to an efficiency of
$70\%$.

\section{Conclusions}\label{sec:conclusions}

In this work we have proposed and analyzed the numerical
performance of parallelization strategies for BoostDMS, an
implementation of a Direct Multisearch method.

Additionally to the obvious parallelization of function evaluation
at the poll step, strategies for parallelizing the function
evaluation at the search step were also considered, stressing the
benefits associated to the evaluation by levels performed in the
sequential version of the code. The value of parallel strategies
for models computation and minimization at the search step was
also assessed. New strategies for the selection of the iterate
point were proposed, based on modifications of the $\Gamma$ spread
metric, and considering the possibility of selecting several
iterate points at a given iteration.

The best strategy, at each iteration selects $q+2$ iterate points
from the current approximation to the Pareto front, corresponding
to the best value for each component of the objective function and
to the two points defining the largest gap in a given component of
the objective function, with cyclic changes between iterations.
Function evaluation is always performed in parallel, respecting
the levels approach of the search step.

This best version allows a remarkable improvement in computational
time, with a significant improvement in purity and a slight
improvement in hypervolume metrics. Conclusions hold both in an
extensive academic test set and in a real engineering application.

\section*{Funding}
Support for authors was provided by national funds through FCT --
Funda\c c\~ao para a Ci\^encia e a Tecnologia I. P., under the
scope of projects PTDC/MAT-APL/28400/2017 and UIDB/MAT/00297/2020,
the later only for second and third authors.

\bibliographystyle{tfs}
\bibliography{bibboostparallel}

\end{document}